\documentclass[11pt]{article}
\pdfoutput=1
\usepackage{amsmath}    
\usepackage{graphicx}   
\usepackage{verbatim}   
\usepackage{color}      
\usepackage{subfigure}  
\usepackage{hyperref}   
\usepackage{bbm}

\usepackage{amssymb,amsmath,amsthm,a4wide}
\newcommand{\Rset}{\mathbb{R}}
\newcommand{\Cset}{\mathbb{C}}
\newcommand{\Nset}{\mathbb{N}}
\newcommand{\primal}{\mathbf{(P)}}
\newcommand{\dual}{\mathbf{(D)}}
\newcommand{\saddleP}{\mathbf{(S)}}

\newcommand{\penaltyT}{R}
\newcommand{\dataMisfit}{S}
\newcommand{\operator}{T}
\newcommand{\X}{X}
\newcommand{\q}{}
\newcommand{\qstar}{}
\newcommand{\Y}{Y}

\newcommand{\CP}{\text{CP}}
\newcommand{\CPBS}{\text{CP-BS}}

\newcommand{\xbar}{\bar{x}}
\newcommand{\ybar}{\bar{y}}

\newcommand{\xhat}{\hat{x}}
\newcommand{\pbar}{\bar{p}}

\newcommand{\FT}{\mathcal{F}}

\newcommand{\obfct}{\mathcal{O}}

\newcommand{\ima}{\text{i}}
\newcommand{\e}{\text{e}}

\newcommand{\torus}{\mathbb{T}}

\newcommand{\BregD}{\mathcal{B}}
\DeclareMathOperator*{\argmin}{argmin}
\DeclareMathOperator*{\argmax}{argmax}

\usepackage{amsthm}
\usepackage{amsmath}
\usepackage{amssymb}
\usepackage{amsfonts}
\usepackage{nicefrac}
\usepackage{txfonts}
\usepackage{mathrsfs}
\usepackage{array}
\usepackage{graphicx}
\usepackage{pstricks-add}
\usepackage{textcomp}

\theoremstyle{definition}
\newtheorem{mydef}{Definition}
\newtheorem{myexp}[mydef]{Example}
\theoremstyle{plain}
\newenvironment{myproof}{\begin{proof}[\large{\textbf{Proof.}}]}{\end{proof}}
\newtheorem{mythm}[mydef]{Theorem}
\newtheorem{mycor}[mydef]{Corollary}
\newtheorem{lemma}[mydef]{Lemma}
\theoremstyle{remark}

\newtheorem{myremark}[mydef]{Remark}
\newtheorem{myalg}{Algorithm}
\title{A Generalization of the Chambolle-Pock Algorithm to Banach Spaces
with Applications to Inverse Problems}
\author{Thorsten Hohage and Carolin Homann}

\begin{document}
\maketitle
\begin{abstract}
\noindent
For a Hilbert space setting Chambolle and Pock introduced an attractive first-order algorithm  
which solves a convex optimization problem and its Fenchel dual simultaneously. 
We present a generalization of this algorithm to Banach spaces.
Moreover, under certain conditions we prove strong convergence 
as well as convergence rates. Due to the generalization the method becomes efficiently
applicable for a wider class of problems.
This fact makes it particularly interesting for solving ill-posed inverse problems 
on Banach spaces by Tikhonov regularization or the iteratively regularized Newton-type method, respectively. 
\end{abstract}

\section{Introduction}
Let \( \X \) and \( \Y \) be real Banach spaces 
and \( \operator: \X \to \Y \) a linear, continuous operator, 
with adjoint \( \operator^*:\Y^* \to \X^*\). 
In this paper we will consider the convex optimization problem
\begin{equation} \label{prbl:primal}
\bar{x} = \text{argmin}_{x \in \X} \left( g(\operator x) + f(x) \right) \qquad \primal 
\end{equation}
as well as its Fenchel dual problem 
\begin{equation}\label{prbl:dual}
 \bar{p} =  \argmax \limits_{p \in Y^*} \left( -f^*(\operator^*p) - g^*(-p) \right), \qquad \dual
\end{equation}
where \( f:\X \to [0, +\infty) \) and \( g: \Y \to [0, +\infty) \) belong to the class 
\( \Gamma(X)\) and \( \Gamma(Y) \) of proper, convex and lower semicontinuous (l.s.c.) functions. 
By \(f^* \in \Gamma(\X^*)\) and \( g^* \in \Gamma(\Y^*) \) we denote their conjugate functions. 
A problem of the form \( \primal \) arises in many applications, such as
image deblurring (e.g. the ROF model \cite{ROFmodel}), sparse signal restoration 
(e.g. the LASSO problem \cite{LASSO}) and inverse problems.
We would like to focus on the last aspect. 
Namely solving a linear ill-posed problem \( \operator x = y \) 
by the Tikhonov-type regularization of the general form 
\begin{equation} \label{eq:Tikhonov_functional}
x_{\alpha} =  \argmin \limits_{x \in X} \dataMisfit(y;\operator x) + \alpha \ \penaltyT(x) ,
\end{equation}
leads for common choices of the data fidelity functional \( \dataMisfit(y; \cdot)\in \Gamma(\Y) \) and
the penalty term \( \penaltyT \in \Gamma(\X) \) to a problem of the form \( \primal \). 
Also for a nonlinear  Fr\'{e}chet differentiable  operator \( T:X \to Y ,\)  
where the solution of the operator equation
\( T x = y \) can be recovered by the iteratively regularized Newton-type method (IRNM, see e.g. \cite{HohageWerner})
\begin{equation} \label{eq:IRNM}
x_{n+1} =  \argmin \limits_{x \in X} S(y; T(x_n) + T^{\prime}[x_n] (x - x_n)) + \alpha_{n} \ R(x), 
\end{equation}  
we obtain a minimization problem of this kind in every iteration step.
In particular if \( \dataMisfit \) or \( \penaltyT \) is 
(up to an exponent) given by a norm  \( \| \cdot \|_Z \) of a Banach space \( Z \) 
it seems to be natural to choose \( \X \) respectively \( \Y \) equal to \( Z \).
These special problems are of interest in the current research, 
see e.g.\cite{KaltenbacherHofmann, SchoepferLouisSchuster, schuster2012regularization} and references therein.\\
Also inverse problems with Poisson data, which occur for example in photonic imaging,
are a topic of certain interest (cf. \cite{HohageWerner, WernerHohage}).
Due to the special form of the appropriate choice of \( S \), 
here proximal algorithms appear to be particularly suitable 
for minimizing the corresponding regularization functional.

If \( \X \) and \( \Y \) are Hilbert spaces
one finds a wide class of first-order proximal algorithms in literature for solving \( \primal \) , 
e.g. FISTA \cite{BeckTeboulle} , ADMM \cite{Boyd_ADMM}, proximal splitting algorithms \cite{CombettesPesquetOverview}.
Chambolle and Pock introduced the following first-order primal-dual algorithm (\cite{ChambollePock}), 
which solves the primal problem \( \primal\) and its dual \( \dual \) simultaneously:
\begin{myalg}[CP] 
For suitable choices of \( \tau_k, \sigma_k > 0,\ \theta_k \in [0,1], \left(x_0,p_0 \right) \in X \times Y, \hat{x}_0 \coloneqq x_0 \), set:
\begin{align}
\label{alg:CP_1}
p_{k+1} &= \left( \sigma_k \, \partial g^* + I \right)^{-1} \left(p_k + \sigma_k \operator \hat{x}_k \right)\\
\label{alg:CP_2}
x_{k+1} &= \left( \tau_k \, \partial f + I \right)^{-1} \left( x_k -\tau_k \operator^*p_{k+1}\right)\\
\label{alg:CP_3}
\hat{x}_{k+1} &= x_{k+1} + \theta_k \left(x_{k+1}-x_k\right)
\end{align}
\end{myalg}

\noindent
Here \( \partial f \) denotes the (set-valued) subdifferential of a function \( f ,\)
which will be defined in section \ref{sec:preliminaries}. 
There exists generalizations of this algorithm in order to solve monotone inclusion 
problems (\cite{Bot_PD_Monotone, Vu_PD_COCOERCIVE}) 
and to the case of nonlinear operators \( T \) (\cite{Valkonen_nonlinearPD}).
Recently, Lorenz and Pock (\cite{LorenzPock14}) proposed a quite
general forward-backward algorithm for monotone inclusion 
problems with \( \CP \) as an special case.

In \cite{ChambollePock} there are three different parameter choice rules given,
for which strong convergence was proven. Two of them base on the assumption
that \( f \) and/or \(g^* \) satisfy a convex property, which enable to prove convergence rates.
In order to speed up the convergence of the algorithm, 
\cite{ChambollePock_precond} and \cite{He_Yuan_preconCP} 
discuss efficient preconditioning techniques. 
Thereby the approach studied in \cite{ChambollePock_precond} can be
seen as a generalization from Hilbert spaces \( \X \) and \( \Y \) 
to spaces of the form \(\Upsilon^{\frac{1}{2}} X \) and \( \Sigma^{-\frac{1}{2}}Y \)
for symmetric, positive definite matrices \( \Upsilon \) and \( \Sigma \), 
where the dual spaces with respect to standard scalar products are given by
\(\Upsilon^{-\frac{1}{2}} X \) and \( \Sigma^{\frac{1}{2}}Y \), respectively.
Motivated by this approach, in this paper we further develop a nonlinear generalization of \( \CP \) 
to reflexive, smooth and convex Banach spaces \( \X \) and \( \Y \), 
where we 
assume \( \X \) to be \( 2\)-convex and \( \Y\) to be \(2\)-smooth. 
For all three variations of \( \CP \), introduced in \cite{ChambollePock}, 
we will prove the same convergence results, including linear convergence
for the case that \( f \) and \(g^* \) satisfy a specific convex property on Banach spaces.
Moreover the generalization provides clear benefits regarding the efficiency and the feasibility:

\noindent 
First of all the essential factor affecting the performance of the \(\CP \)-algorithm is the efficiency 
of calculating the (well-defined, single-valued) resolvents \( \left( \sigma \, \partial g^* + I \right)^{-1} \) and 
\( \left( \tau \, \partial f + I \right)^{-1} \) 
(cf. \cite{Villa_accelerated} addressing problem of non exact resolvents in forward-backward algorithms). 
By the generalization of \(\CP \) and of the resolvents, inter alia, 
we obtain closed forms of these  operators for a wider class of functions \( f \) and \( g\).
Furthermore, there exists a more general set of functions that fulfill the generalized convex property 
on which the accelerated variations of \( \CP \) are based on.
Moreover, in numerical experiments we obtained faster convergence
for appropriate choices of \( X\) and \( Y \).

The paper is organized as follows: In the next section we give necessary definitions and results of convex analysis
and optimization on Banach spaces. In section \ref{sec:alg_convergence_re} we present a generalization
of \(\CP\) to Banach spaces, and prove convergence results for special parameter choice rules. 
The generalized resolvents, which are included in the algorithm, are the topic of section \ref{sec:general_resolvents}.
In order to illustrate the numerical performance of the proposed method we apply it in section \ref{sec:numerical_exam}
to some inverse problems. In particular we consider a problem with sparse solution and a special phase retrieval problem, 
given as a nonlinear inverse problem with Poisson data.

\section{Preliminaries} \label{sec:preliminaries}
The following definitions and results from convex optimization and geometry of Banach spaces 
can be found e.g. in \cite{BarbuPrecupanu, Cioranescu}.

For 
a Banach space \( Z \) let  \( Z^* \) denote its topological dual space. 
In analogy to the inner product on Hilbert spaces, we write
\( \left<z,z^* \right>_Z = z^*(z) \) for \( z \in Z \) and \( z^* \in Z^* \).
Moreover, for a 
function \( h \in \Gamma(Z) \) on a Banach space \( Z \) let 
\( \partial h: Z \rightrightarrows Z^*\),  
\(z \mapsto \left\{ z^* \in  Z^* \, | \ \forall u \in Z\;\left< u-z, z^* \right>_Z \leq h(u) -h(z) \right\} \)
denote the subdifferential of \( h\). 
Then \( \xbar\) is the unique minimizer  of \( h \in \Gamma(\X) \) if and only if
\( 0 \in \partial h(\xbar) \).
Moreover, under certain conditions \( \xbar \in \X \) is a solution to the primal problem \( \primal \) 
and  \( -\pbar \in \Y^* \) is a solution to the dual problem \( \dual \) 
if and only if the optimality conditions (see e.g.  \cite[Section 2.8]{Zalinescu_cvxAna})
\begin{equation} \label{eq:extremal_relation}
 -\operator^* \pbar \in \partial f(\xbar), \quad \operator \xbar \in \partial g^{*}(\pbar)
\end{equation}
hold. Another equivalent formulation is that the pair \( \left( \xbar, -\pbar \right) \) 
solves the saddle-point problem
\[\min \limits_{x \in \X} \max \limits_{p \in \Y^*} -\left< \operator x, p \right>_{\Y} + f(x) - g^*(-p) \quad \saddleP.\]
Rewriting lines \eqref{alg:CP_1} and \eqref{alg:CP_2} as
\[ \sigma_k\, T \hat{x_k} \in  \sigma_k\,  \partial g^* + p_{k+1} - p_{k}, 
\quad - \tau_k \operator^*p_{k+1} \in \tau_k \, \partial f(x_{k+1}) + x_{k+1} - x_{k}  \]
we can interpret the \( \CP\)-algorithm as a fixed point iteration with   
an over-relaxation step in line \eqref{alg:CP_3}. 
The objective value can be expressed by the \emph{partial primal-dual gap}: 
\[ \mathcal{G}_{B_1 \times B_2}(x,p) \coloneqq 
\max \limits_{p^{\prime} \in B_2} \left(\left< \operator x, p^{\prime} \right>_{\Y} - g^*(p^{\prime}) + f(x) \right)
- \min \limits_{x^{\prime} \in B_1} \left( \left< \operator x^{\prime}, p \right>_{\Y} - g^*(p) + f(x^{\prime}) \right),\]
on a bounded subset \(B_1 \times B_2 \subset \X \times \Y^*\) . 

\paragraph*{Resolvents}
On a Hilbert space \( Z \) the operator \( \left( \tau \partial h + I \right) \) 
is bijective for any \( h \in \Gamma(Z) \) and any \( \tau > 0 \),  
i.e. the resolvent 
\( \left( \tau \partial h + I \right)^{-1}: Z^* = Z \to Z\) of \( h \) 
is well-defined and single-valued. More generally Rockafellar proved (\cite[Proposition 1]{Rockafellar_resolvent}) 
that on any reflexive Banach space \( Z \) 
the function \(\left( \tau \partial h + J_Z \right)^{-1}: Z^* \to Z,\) 
where \( J_Z = \partial \Phi \) is the subdifferential of \( \Phi(x) = \frac{1}{2} \|x\|_Z^2\),
is well-defined and single-valued, as well. 
Furthermore, 
as \( z = \left( \tau \, \partial h + J_Z \right)^{-1} (u) \) is the unique solution of
\[ 0 \in \partial\, \tau h(z) + \partial  \Phi(z) - u \]
this generalized resolvent can be rewritten as follows:
\begin{equation} \label{eq:Gresolvent_argmin}
\left( \tau \, \partial h + J_Z \right)^{-1} (u) 
=  \argmin_{z \in Z} \left( \tau \, h (z) - \left<z, u \right>_Z +  \frac{1}{2} \| z \|_Z^2 \right) .
\end{equation}

\paragraph*{Regularity of Banach spaces}
We make some assumptions on the regularity of the Banach spaces \( \X \) and \( \Y \).
A Banach space \( Z \) is said to be \( r-\)convex with \( r > 1\) if there exists a constant \( C > 0 \),
such that the \emph{modulus of convexity} \( \delta_Z:[0,2] \to [0,1],\)
\[ \delta_Z (\epsilon ) \coloneqq \inf \left\{ \left. 1- \left\| \frac{1}{2} \left(x+u \right) \right\|_Z \, 
\right| \ \left\|x\right\|_Z = \left\|u \right\|_Z = 1, \left\|x-u\right\|_Z \leq \epsilon \right\} \]
satisfies 
\[ \delta_Z(\epsilon) \geq C\, \epsilon^r, \quad \epsilon \in [0,2].\]
We call \( Z \) \( r-\)smooth, if the \emph{modulus of smoothness} \(\rho_Z: [0, \infty] \to \Rset \) 
\[ \rho_X (\tau ) \coloneqq \frac{1}{2} \sup \left\{ \left.  \left\| x+u  \right\|_Z -\left\| x-u  \right\|_Z -2
\, \right| \ \left\|x\right\|_Z = 1, \left\|u\right\|_Z \leq \tau \right\} \]
fulfills the inequality
\(\rho_X (\tau ) \leq G_Z\, \tau^r \)
for any \(\tau \in [0, \infty]\) and some constant \( G_Z > 0 \).
In the following we will assume both  \( \X \) and the dual space \(Y^* \) to be reflexive, smooth and 2-convex Banach spaces. 
Because of the Lindenstrauss duality formula the second statement is equivalent to the condition that \( Y \) is a reflexive, 
convex and 2-smooth Banach space. \\

\paragraph*{Duality mapping}
For \( q \in (1, \infty) \) let us introduce the \emph{duality mapping}
\[ J_{q,Z}: Z \rightrightarrows Z^* , J_{q,Z}(x) 
\coloneqq \left\{ x^* \in Z^* \,|\ \left<x,x^* \right>_Z = \left\| x \right\|_Z \  \left\| x^* \right\|_{Z^*},  
\left\| x^* \right\|_{Z^*} = \left\| x \right\|_Z^{q-1}  \right\} \]  
with respect to the weight function \( \phi(t) =  t^{q-1}. \)
If \( Z \) is smooth, \( J_{q,Z} \) is single-valued. 
If in addition \( Z \) is \(2\)-convex and reflexive, 
\( J_{q,Z} \) is bijective with inverse \( J_{q^*, Z^*}: Z^* \to Z^{**} = Z, \) 
where \( q^* \in (1, \infty),\)  denotes the conjugate exponent of \( q \), i.e. \( \frac{1}{q} + \frac{1}{q^*} = 1 \).
By the theorem of Asplund (see \cite{Asplund}) \( J_{q,Z} \) can be also defined
as the subdifferential \( \partial \Phi_q \) of  \( \Phi_q(x) \coloneqq \frac{1}{q} \|x\|_Z^q \).
Thus, for the case \( q = 2 \) the so-called normalized duality mapping \( J_{2,Z} \) 
coincides with the function \( J_Z \), we introduced in the previous section.
Note that the duality mapping is in general nonlinear.

\paragraph*{Bregman distance}
Instead of for the functional \( (u,v) \mapsto \|u-v\|_X^2\) 
we will prove our convergence results with respect to the \emph{Bregman distance}
\[ \BregD_{Z}(u,x) \coloneqq 
\frac{1}{2} \left\| u \right\|_Z^{2} - \frac{1}{2} \left\| x \right\|_Z^{2} - \left< u -x, J_{Z} (x) \right>_Z  \] 
with gauge function \( \Phi_{2} \). 
Note that \(  \BregD_{Z} \) is not a metric, as symmetry is not fulfilled. 
Nevertheless, a kind of symmetry with respect to the duals holds true:
\begin{align}
\label{eq:Bregman_dualSymmetrie}
 \BregD_{Z^*} \left(J_{Z} \left(v \right), J_{Z} \left(x \right) \right) 
 = \BregD_{Z} \left(x, v \right).
\end{align}
Moreover, the Bregman distance \(\BregD_{Z}\) satisfies the following identity
\begin{align}
\label{eq:triag_bregman_1}
\begin{split}
 \BregD_{Z}\left(u,x\right) + \BregD_{Z} \left(v,u \right) 
 &= \frac{1}{2} \left\| v \right\|_Z^{2} - \frac{1}{2} \left\| x \right\|_Z^{2} - 
 \left< u -x, J_{Z} (x) \right>_Z - \left< v -u, J_{Z} (u) \right>_Z \\
 &=  \BregD_{Z} \left(v,x\right) + \left< v -u, J_{Z} (x) - J_{Z} (u) \right>_Z, \quad  x,u,v \in Z.
 \end{split}
\end{align}
From \eqref{eq:Gresolvent_argmin} we observe that 
\[ \left( \tau \, \partial h + J_Z \right)^{-1} (u) 
= \argmin_{z \in Z} \left( \tau \, h (z) + \BregD_{Z}(z, J_{Z^*}(u) \right) \]
is a resolvent with respect to the Bregman distance.\\

The assumption \( \X \) and \(\Y^* \) to be reflexive, smooth and 2-convex Banach spaces 
provide the following helpful inequalities (see e.g. \cite{MassSchoepferSchuster}):
There exist positive constants \( C_{\X}\) and \( C_{\Y^*} \), such that:
\begin{equation} \label{ineq:Bregman_pNorm}
\mathcal{B}_{\X} \left(x, u \right) \geq \frac{C_{X}}{2} \|x - u \|_{\X}^2\ \quad \forall x,u\in \X, 
\qquad \text{and} \qquad 
\mathcal{B}_{\Y^*}\left( y^*, p \right)  \geq \frac{C_{\Y^*}}{2} \|y^* - p \|_{\Y^*}^{2} \quad \forall y^*,p \in \Y^*.
\end{equation}

\begin{myexp}
 Considering the proof of the last inequalities \eqref{ineq:Bregman_pNorm}, we find that the constant \( C_{\X} >0 \) 
comes from a consequence of the Xu-Roach inequalities (\cite{XuRoach}):  
 \begin{equation} \label{ineq:estimate_C_X}
\frac{1}{2} \| x -u \|_{\X}^2 \geq \frac{1}{2} \|x\|_{\X}^2 -  \left< J_{\X}x,u \right>_{\X} + \frac{C_{\X}}{2} \|u \|_{\X}^2,\quad u,x \in \X.
\end{equation}
For example for \( \X = l^r \) with \( r \in (1,2] \) this estimate holds for \( C_{\X} < r-1\) as it is shown in \cite{Xu_Lp_ineq}.
\end{myexp}

\begin{lemma}\label{cor:Young_BS}
 For any \( x, u \in \X \), \( y^*, p \in \Y^* \) and any positive constant \( \alpha \), we have
 \begin{equation}\label{ineq:Young_BS}
  \left|\left< \operator(x-u), p-y^*\right>_{\Y} \right|
  \leq \|\operator\| \left(\frac{\alpha \, \min \left\{\BregD_{\X}\left(x,u\right), \BregD_{\X}\left(u,x\right) \right\}}{C_{\X}} 
  +  \frac{\min \left\{ \BregD_{\Y^*}\left(p, y^* \right), \BregD_{\Y^*}\left( y^*, p \right) \right\}}{\alpha\, C_{\Y^*}}\right),
 \end{equation}
 where \( \|\operator \| = \max \left\{ \|\operator x \|_{\Y}    \, | \ x \in \X,  \|x\|_{\X} = 1 \right\}  \) denotes the operator norm.
\end{lemma}
\begin{myproof}
Applying Cauchy-Schwarz's inequality as well as
the special case of Young's inequality:
\begin{equation} \label{ineq:Young}  
ab \leq \frac{\alpha\, a^2}{2} + \frac{b^{2}}{2\, \alpha}, \quad a,b \geq 0
\end{equation} 
with \( a  \coloneqq \| x_{k} - x_{k-1} \|_{\X},\) and \hbox{\(b \coloneqq \|p_{k+1} - p_k \|_{\Y^*} \in \Rset \)}
leads to
\begin{align*} 
\left|\left< \operator(x-u), p-y^*\right>_{\Y} \right| \leq \| \operator \|\, \| x-u \|_{\X} \, \|p-y^* \|_{\Y^*} 
\leq \| \operator \|\left(\alpha\, \frac{\| x-u \|_{\X}^2}{2} + \frac{\|p-y^*\|_{\Y^*}^2}{2\, \alpha}\right).
\end{align*}
Now, the inequalities \eqref{ineq:Bregman_pNorm} gives the assertion.
\end{myproof}

\section{Algorithms and convergence results} \label{sec:alg_convergence_re}
\begin{myalg}[CP-BS] \label{alg:CP_BS1}
For \(  \left(\tau_k, \sigma_k\right)_{k \in \Nset}\subseteq (0, \infty),
\ \theta \in [0,1], \left(x_0,p_0 \right) \in \X \times \Y^*, 
\hat{x}_0 \coloneqq x_0 \), set: 
\begin{align}
\label{alg:CP_BS_1}
p_{k+1}  &\coloneqq \left(  \sigma_k \,  \partial g^* +  J_{Y^*} \right)^{-1}
\left( J_{Y^*} \left(  p_k \right) +  \sigma_k\, \operator \xhat_k \right) \\
\label{alg:CP_BS_2}
x_{k+1}  &\coloneqq \left(\tau_k\,  \partial  f + J_{X} \right)^{-1}  
\left( J_{X} \left( x_k \right) - \tau_k \,   \operator^* p_{k+1} \right) \\
\label{alg:CP_BS_3}
\xhat_{k+1} &\coloneqq x_{k+1} + \theta_k \left(x_{k+1} -  x_{k} \right)
\end{align} 
\end{myalg}
Let us assume that there is a solution \( \left( \xbar, -\pbar \right) \) to the saddle-point problem \( \saddleP \).
In analogy to \cite{ChambollePock} we like to bound the distance 
of one element of the sequence \( \left(x_k, p_{k} \right)_{k \in \Nset} \) to a solution of \( \saddleP \). 
For the given general Banach space case 
we measure this misfit by Bregman distances and define for an arbitrary point 
\(\left(x, y^* \right) \in X \times Y^* \)  
\[ \triangle_k(x, y^*) \coloneqq \frac{\BregD_{Y^*}(y^*,p_k)}{\sigma_k} +  \frac{\BregD_{X}(x, x_k)}{\tau_k} .\]
\begin{mythm}\label{thm:CP_BS1}
We choose constant parameters \(\sigma_k = \sigma,\) \( \tau_k = \tau \) and \( \theta_k = 1\) 
for some \( \sigma, \tau \) with
\( \sqrt{\sigma  \, \tau} \| \operator \| < \text{min}\, \left\{C_X, C_{Y^*} \right\}  \), 
where \(C_X, C_{Y^*}\) are given by \eqref{ineq:Bregman_pNorm}.
Then for algorithm \ref{alg:CP_BS1} the following assertions hold true:
\begin{itemize}
 \item The sequence \( \left(x_k, p_k \right)_{k \in \Nset} \) remains bounded. 
 More precisely there exists a constant 
 \[ C < \left(1- \frac{\| \operator \|^{2} \, \sigma\, \tau}{C_X\, C_{Y^*}} \right)^{-1},\] 
 such that for any \( N \in \Nset \)
 \begin{equation} \label{ineq:xk_pk_bounded_wrt_B}
 \triangle_N(\xbar, \pbar)  \leq  C \ \triangle_0(\xbar, \pbar). 
\end{equation}
 \item The restricted primal-dual gab \(\mathcal{G}_{B_1 \times B_2} \left(x^N, p^N \right) \) at the mean values 
 \( x^N \coloneqq \frac{1}{N} \sum_{k = 1}^N x_k \in X\) and 
 \( y^N \coloneqq \frac{1}{N} \sum_{k = 1}^N y_k\in Y^* \) is bounded by
 \[
 D(B_1, B_2) \coloneqq \frac{1}{N} \sup \limits_{(x, y^*) \in B_1 \times B_2} 
 \triangle_0 \left(y^*, x\right) 
 \]
  for any bounded set \(B_1 \times B_2 \in X \times Y^* \). 
  Moreover, for every weak cluster point \( \left( \tilde{x}, \tilde{p} \right)\) 
  of the sequence \( \left( x^N, p^N \right)_{N \in \Nset} \), 
  \( \left( \tilde{x}, -\tilde{p} \right)\) solves the saddle-point problem \( \saddleP \).
  \item If we further assume the Banach spaces \( X \) and \( Y\) to be finite dimensional, 
  then there exists a solution \( \left( \xbar, -\pbar \right) \)
  to the saddle-point problem \( \saddleP \) such that the sequence \( \left(x_k, p_k \right) \) 
  converges strongly to \( \left( \xbar, \pbar \right)\).
\end{itemize}
\end{mythm}

\begin{myproof}
 
Using \eqref{eq:triag_bregman_1} 
the misfit functional \( \triangle_k(x, y^*) \) for  some \( \left( x, y^* \right) \in X \times Y^* \) reads
\begin{align*} 
&\triangle_k(x, y^*)
= \frac{\BregD_{Y}\left(J_{Y^*} \left(p_{k} \right),J_{Y^*} \left( y^* \right)\right)}{ \sigma} 
+ \frac{\BregD_{X^*} \left(J_{X} \left(x_{k} \right), J_{X} \left(x \right) \right)}{\tau} \\
=& \frac{\BregD_{Y}\left(J_{Y^*} \left(p_{k+1} \right), J_{Y^*} \left(y^* \right)\right)}{ \sigma} - 
\left< \frac{J_{Y^*} \left(p_{k} \right) - J_{Y^*} \left(p_{k+1} \right) }{\sigma}, y^* - p_{k+1} \right>_{Y} 
+ \frac{\BregD_{Y} \left( J_{Y^*} \left(p_{k} \right), J_{Y^*} \left(p_{k+1} \right) \right)}{ \sigma} \\
&+  \frac{\BregD_{X^*} \left(J_{q,X} \left(x_{k+1} \right), J_{X} \left(x \right) \right)}{\tau} 
- \left< x - x_{k+1},  \frac{J_{X} \left(x_{k} \right) - J_{X} \left(x_{k+1} \right) }{\tau} \right>_{X} 
+ \frac{\BregD_{X^*} \left(J_{X} \left(x_{k} \right), J_{X} \left(x_{k+1} \right) \right)}{\tau}. 
\end{align*}
The iteration formulas \eqref{alg:CP_BS_1} and \eqref{alg:CP_BS_2} imply that
\[ \frac{1}{\sigma} \left(J_{Y^*}\left(p_k \right)- J_{Y^*}\left(p_{k+1} \right) \right)+ \operator \xhat_k  \in \partial g^*(p_{k+1})
\quad \text{and} \quad
 \frac{1}{\tau} \left(J_{X}( x_k ) - J_{X} \left(x_{k+1} \right) \right) - \operator^* p_{k+1}  \in f(x_{k+1}). \]
So by the definition of the subdifferential we obtain: 
\begin{align}
\label{eq:def_subd_for_p_k+1}
 g^* \left(y^* \right) - g^* \left(p_{k+1} \right) &\geq  
 \left< \frac{ J_{Y^*}\left(p_k \right)- J_{Y^*} \left(p_{k+1} \right)}{\sigma}, y^* - p_{k+1} \right>_Y + 
\left< \operator \xhat_k, y^*-p_{k+1} \right>_{Y}  \\ 
\label{eq:def_subd_for_x_k+1}
f \left( x\right) - f \left( x_{k+1}\right) &\geq
\left< x-x_{k+1}, \frac{ J_{X}( x_k ) - J_{X} \left(x_{k+1} \right) }{\tau} \right>_{X}
-\left< \operator \left(x-x_{k+1}\right),  p_{k+1}  \right>_X ,
\end{align}
Using \eqref{eq:Bregman_dualSymmetrie} this yields
\begin{align} 
\label{ineq:CP_12_weighted_error}
\triangle_k(x, y^*)
\geq &\  g^* \left(p_{k+1} \right) -g^* \left(y^* \right)  - \left< \operator \xhat_k, p_{k+1} - y^* \right>_{Y}
+ f \left( x_{k+1}\right) - f \left( x\right) -\left< \operator \left(x_{k+1} - x\right),  - p_{k+1}  \right>_X \\
&+ \frac{\BregD_{Y^*}\left(y^*, p_{k+1} \right)}{ \sigma} + \frac{\BregD_{Y^*} \left( p_{k+1}, p_{k}  \right)}{ \sigma} 
+  \frac{\BregD_{X} \left( x, x_{k+1} \right)}{\tau} + \frac{\BregD_{X} \left(x_{k+1}, x_{k} \right)}{\tau} \\
\ & + \left< \operator x_{k+1}, p_{k+1} - y^* \right>_Y  - \left< \operator \left(x_{k+1} -x \right), p_{k+1} \right>_Y 
+ \left< \operator x_{k+1}, y^* \right>_Y - \left< \operator x, p_{k+1} \right>_Y \\
\label{ineq:PDGab_in_CP_12_2}
= & \left[\left< \operator x_{k+1}, y^* \right>_Y - g^*( y^*) + f(x_{k+1}) \right] 
- \left[\left< \operator x, p_{k+1} \right>_Y - g^*(p_{k+1}) + f(x) \right] \\
& + \triangle_{k+1}(x, y^*)
+ \frac{\BregD_{Y^*} \left( p_{k+1}, p_{k}  \right)}{ \sigma} 
+ \frac{\BregD_{X} \left(x_{k+1}, x_{k} \right)}{\tau}
\label{ineq:CP_12_r}
+ \left< \operator (x_{k+1} - \xhat_k ), p_{k+1} - y^* \right>_Y. 
\end{align}
In order to estimate the last summand in \eqref{ineq:CP_12_r}, we insert \eqref{alg:CP_BS_3} 
with \( \theta_k = 1 \),
and apply Lemma \ref{cor:Young_BS}  
with
\( \alpha \coloneqq \left(\frac{\sigma}{\tau}\right)^{\frac{1}{2}} > 0 \):
\begin{align*}
 &\left< \operator \left((x_{k+1} - x_{k}) - (x_{k} - x_{k-1}) \right), p_{k+1} - y^* \right>_Y \\
 &= \left< \operator \left(x_{k+1} - x_{k}\right), p_{k+1} - y^* \right>_Y - \left< \operator \left(x_{k} - x_{k-1} \right), p_{k} - y^* \right>_Y 
 - \left< \operator \left( x_{k} - x_{k-1} \right), p_{k+1} - p_k \right>_Y \\
 &\geq \left< \operator \left(x_{k+1} - x_{k}\right), p_{k+1} - y^* \right>_Y - \left< \operator \left(x_{k} - x_{k-1} \right), p_{k} - y^* \right>_Y \\
 &- \frac{\|\operator \|\, \sigma^{\frac{1}{2}}\, \tau^{\frac{1}{2}}}{C_X} \ \frac{\BregD_{\X} \left( x_{k},x_{k-1} \right)}{\tau} 
  - \frac{\|\operator \|\, \sigma^{\frac{1}{2}}\, \tau^{\frac{1}{2}}}{C_{\Y^*}}\ \frac{ \BregD_{\Y^*} \left(p_{k+1}, p_k \right)}{\sigma}.
\end{align*}
Thus, we conclude that
\begin{align*} 
\triangle_k(x, y^*) 
\geq & \left[\left< \operator x_{k+1}, y^* \right>_Y - g^*( y^*) + f(x_{k+1}) \right] - \left[\left< \operator x, p_{k+1} \right>_Y - g^*(p_{k+1}) + f(x) \right] \\
&+ \triangle_{k+1} (x, y^*) 
+ \left(1- \frac{\| \operator \|\, \sigma^{\frac{1}{2}}\, \tau^{\frac{1}{2}}}{C_{Y^*}} \right)\frac{\BregD_{Y^*} \left( p_{k+1}, p_{k}  \right)}{ \sigma} 
 \ - \frac{\| \operator \|\, \sigma^{\frac{1}{2}}\, \tau^{\frac{1}{2}}}{C_{X}} \frac{\BregD_{X} \left(x_{k}, x_{k-1} \right)}{\tau}\\
\ &+ \frac{\BregD_{X} \left(x_{k+1}, x_{k} \right)}{\tau} 
+ \left< \operator \left(x_{k+1} - x_{k}\right), p_{k+1} - y^* \right>_Y - \left< \operator \left(x_{k} - x_{k-1} \right), p_{k} - y^* \right>_Y .
\end{align*}
Summing from \( k = 0 \) to \(N-1\) leads to
\begin{align*}
\triangle_0(x, y^*)  &  + | \left< \operator \left(x_{N} - x_{N-1} \right), p_{N} - y^* \right>_Y |\\
&\geq \sum_{k = 0}^{N} \left[\left< \operator x_{k+1}, y^* \right>_Y - g^*( y^*) + f(x_{k+1}) \right] 
- \left[\left< \operator x, p_{k+1} \right>_Y - g^*(p_{k+1}) + f(x) \right] \\
&+ \triangle_{N} (x, y^*) 
+ \left(1- \frac{\| \operator \|\, \sigma^{\frac{1}{2}}\, \tau^{\frac{1}{2}}}{C_{Y^*}} \right) \sum_{k = 1}^{N} \frac{\BregD_{Y^*} \left( p_{k}, p_{k-1}  \right)}{ \sigma} 
+ \frac{\BregD_{X} \left(x_{N}, x_{N-1} \right)}{\tau} \\
&+\left(1- \frac{\| \operator \|\, \sigma^{\frac{1}{2}}\, \tau^{\frac{1}{2}}}{C_{X}} \right) \sum_{k=1}^{N-1} \frac{\BregD_{X} \left(x_{k}, x_{k-1} \right)}{\tau}.
\end{align*}
Now, applying again Lemma \ref{cor:Young_BS}
with \( \alpha = \frac{C_x}{\|\operator\|\, \tau} \) yields
\begin{equation}\label{ineq:apply_Young_ineq_sum}
| \left< \operator \left(x_{N} - x_{N-1} \right), p_{N} - y^* \right>_Y | 
\leq \frac{\BregD_{\X} \left(x_{N}, x_{N-1} \right)}{\tau}
+\frac{\| \operator \|^{2} \, \sigma\, \tau}{C_X\, C_{Y^*}} \frac{\BregD_{Y^*} \left( y^*, p_{N}  \right)}{ \sigma}, 
\end{equation}
so that we deduce 
\begin{align} \label{eq:sum_pg_gap}
\begin{split}
\triangle_0(x, y^*)
&\geq \sum_{k = 0}^{N} \left[ \left< \operator x_{k+1}, y^* \right>_Y - g^*( y^*) + f(x_{k+1}) \right] 
- \left[\left< \operator x, p_{k+1} \right>_Y - g^*(p_{k+1}) + f(x) \right] \\
&+ \left(1- \frac{\| \operator \|^{2} \, \sigma\, \tau}{C_X \, C_{Y^*}} \right)  \frac{\BregD_{Y^*}\left(y^*, p_{N} \right)}{ \sigma} 
+ \left(1- \frac{\| \operator \|\, \sigma^{\frac{1}{2}}\, \tau^{\frac{1}{2}}}{C_{\Y^*}} \right) 
\sum_{k = 1}^{N} \frac{\BregD_{\Y^*} \left( p_{k}, p_{k-1}  \right)}{ \sigma} \\
&+  \frac{\BregD_{\X} \left( x, x_{N} \right)}{\tau} 
+\left(1- \frac{\| \operator \|\, \sigma^{\frac{1}{2}}\, \tau^{\frac{1}{2}}}{C_{\X}} \right) 
\sum_{k=1}^{N-1} \frac{\BregD_{\X} \left(x_{k}, x_{k-1} \right)}{\tau} .
\end{split}
\end{align}
Here, because of the choice 
\( \sigma^{\frac{1}{2}}\, \tau^{\frac{1}{2}} < \frac{\text{min } \left( C_{\X}, C_{\Y^*} \right)}{\|\operator\|} \), 
we obtain only positive coefficients. 
Moreover, for \( \left(x, y^* \right) = \left( \xbar, \pbar \right), \) 
where \( \left( \xbar, -\pbar \right) \) solves the saddle point problem \( \saddleP \),
we have  \( -\operator ^* \pbar \in \partial  f (\xbar) \) and \( \operator \, \xbar \in \partial g^* (\pbar) \),
such that every summand in the first line of \eqref{eq:sum_pg_gap} is non negative as well:
\begin{align} \label{ineq:PDGab_pos}
\begin{split}
 &\left[ -\left< \operator x_{k+1}, -\pbar \right>_Y - g^*( y^*) + f(x_{k+1}) \right] 
 - \left[\left< \operator \xbar, p_{k+1} \right>_Y - g^*(p_{k+1}) + f(\xbar) \right] \\
 &= f(x_{k+1}) - f(\xbar) - \left< x_{k+1} - \xbar, -\operator^* \pbar \right>_X 
+ g^*(p_{k+1}) - g^*( \pbar) - \left< \operator \xbar , p_{k+1} - \pbar \right>_{\Y} 
\geq 0.
\end{split}
\end{align}
This proves the first assertion.
The second follows directly along the lines 
of the corresponding proof in \cite{ChambollePock}, p. 124, 
where we use
\begin{align*} 
\mathcal{G}_{B_1 \times B_2} \left(x^N, p^N \right) 
&= \sup \limits_{(x, y^*) \in B_1 \times B_2  } 
\left[ \left<\operator x_N, y^* \right>_Y - g^*( y^*) + f(x_N) \right] 
- \left[\left< \operator x, p_N \right>_Y - g^*(p_N) + f(x) \right] \\
&\leq \frac{1}{N} \sup \limits_{(x, y^*) \in B_1 \times B_2  } \triangle_0(x, y^*) \end{align*}
instead of (16).
For the last assertion, which needs the assumption that \( \X \) and \( Y \) are finite dimensional, 
we apply again the same arguments as in \cite{ChambollePock}, p. 124, to \eqref{eq:sum_pg_gap} 
from which we obtain
\[ \lim \limits_{k \to \infty} \triangle_k (\xbar, \pbar ) = 0 \]
for a solution \( \left( \xbar, -\pbar \right) \) to the saddle-point problem \( \saddleP \).
This completes the proof.
\end{myproof}

\begin{myremark}
 This generalization covers also the preconditioned version of \( \CP \) proposed in \cite{ChambollePock_precond}: 
There \( \X \) and \( \Y \) are Banach spaces of the form 
 \(\X =\Upsilon^{\frac{1}{2}} H_{\X}\) with \( \|x\|_{\X}= \| \Upsilon^{-\frac{1}{2}} x \|_{H_{\X}} \) 
 and \( \Y =\Sigma^{-\frac{1}{2}} H_{\Y}\) with \( \|y\|_{\Y}= \| \Sigma^{\frac{1}{2}} x \|_{H_{\Y}} \)
for Hilbert spaces \( H_{\X}, H_{\Y} \) and symmetric, positive definite matrices \( \Upsilon \) and \( \Sigma \).
Considering the dual spaces \(\X^* = \Upsilon^{-\frac{1}{2}} H_{\X}  \) and \(\Y^* = \Sigma^{\frac{1}{2}} H_{\Y} \) 
with respect to the scalar product on the corresponding Hilbert spaces, 
the duality mappings read as
\[J_{\X}(x) = \Upsilon^{-1} x,\quad  J_{\Y*}(y) = \Sigma^{-1} \, y .\]
Due to their linearity line \eqref{alg:CP_BS_1} and \eqref{alg:CP_BS_2} take the form of update rule (4)
in \cite{ChambollePock_precond}:
\[p_{k+1} = \left(\sigma_k \Sigma\, \partial g^* + I \right)^{-1} \left(p_k + \sigma_k\, \Sigma \operator \xhat_k \right),
\quad x_{k+1} = \left(\tau_k \Upsilon\, \partial f + I \right)^{-1} \left(x_k - \tau_k\, \Upsilon \operator^* p_{k+1} \right).\]
\end{myremark}

In order to generalize also the accelerated forms of the \( \CP \)-algorithm, 
which base on the assumption that \( f \) is strongly convex,
to Banach spaces, we need a similar property of \( f\). 
More precisely, in \cite{ChambollePock} the following consequence 
of \( f\) being strongly convex with modulus \( \gamma > 0 \) is used:
\begin{align*} 
f(u) - f(x)  
\geq 2\, f\left( \frac{x+u}{2} \right) - 2 f(x) 
+ \frac{\gamma}{2} \|x-u \|_X^2 
\geq  \left< u-x, x^* \right>_{\X} + \frac{\gamma}{2} \|x-u \|_X^2, \quad u,x \in \X, x^*\in \partial f(x).
\end{align*}  
Accordingly, we assume that there exists a constant \(\gamma > 0\) such that \( f \) 
satisfies for any \(x,\, u \in  \X \) and \( x^* \in \partial f(x) \) the following inequality
\begin{equation} \label{ineq:Def_gereneralization_uniformlyCVX}
f\left( u \right)- f(x) \geq \left<u-x, x^* \right>_X + \gamma \ \BregD_{\X}(u, x).
\end{equation}
With this definition we can formulate the next convergence result. 
The case that \eqref{ineq:Def_gereneralization_uniformlyCVX} 
holds for  \(g^* \) instead of \( f \) follows analogously.

\begin{mythm}\label{thm:convergence_CP_BS_2}
Assume that f satisfies \eqref{ineq:Def_gereneralization_uniformlyCVX} for some \( \gamma > 0\) 
and choose the parameters 
\( \left( \sigma_k, \tau_k \right)_{k\in \Nset},\) \(\left( \theta_k \right)_{k\in \Nset} \)
in algorithm \ref{alg:CP_BS1} as follows:
\begin{itemize}
 \item \( \sigma_0 \tau_0 \|T\|^2 \leq \min \left\{C_X, C_{Y^*} \right\} \)
 \item \(\theta_k \coloneqq \left(1+ \gamma\, \tau_{k} \right)^{-\frac{1}{2}}, 
 \quad \tau_{k+1} \coloneqq \theta_{k} \, \tau_k, 
 \quad \sigma_{k+1} \coloneqq \theta_k^{-1} \, \sigma_k \)
\end{itemize}
Then the sequence \( \left(x_k, p_k \right)_{k \in \Nset} \) 
we receive from the algorithm has the following error bound:
For any \( \epsilon > 0 \) there exists an \( N_0 \in \Nset\) such that
\begin{equation}\label{eq:error_bound_CPBS2}
\BregD_{X} \left( \xbar, x_{N} \right)
\leq \frac{4 + 4\,\epsilon}{N^2} 
\left(\frac{\BregD_{X} \left( \xbar, x_{0} \right)}{\gamma^{2} \, \tau_0^2 }
+\frac{\BregD_{Y^*}\left(\pbar, p_{0} \right)}{ \gamma^{2} \, \sigma_0\,\tau_0 }\right) ,
\end{equation}
for all \( N \geq N_0\).
\end{mythm}

\begin{myproof}
We go back to the estimate \eqref{ineq:CP_12_weighted_error}-\eqref{ineq:CP_12_r} 
where we set \( \left(x, y^* \right) \coloneqq \left( \xbar, \pbar \right) \) 
for a solution \( \left( \xbar, -\pbar \right) \) to the saddle-point problem \( \saddleP \).
Assumption \eqref{ineq:Def_gereneralization_uniformlyCVX} applied to 
\( x^* = \nicefrac{1}{\tau_k} \left( J_{\X}(x_k) - J_{X} \left( x_{k+1} \right) \right) 
- \operator^*\, p_{k+1} \in \partial f \left(x_{k+1}\right) \)
gives:
\begin{equation} \label{ineq:subd_fmidcvx}
 f \left( \xbar \right) - f \left( x_{k+1} \right) 
\geq \left< \xbar - x_{k+1}, \frac{J_{\X}(x_k) - J_{X} \left( x_{k+1} \right)}{\tau_k} \right>_{\X} 
- \left< \operator \left(\xbar - x_{k+1} \right),  p_{k+1}  \right>_{\X} 
+  \gamma\, \BregD_{\X}\left(\xbar, x_{k+1} \right).
\end{equation}
Thus replacing \eqref{eq:def_subd_for_x_k+1} by \eqref{ineq:subd_fmidcvx}
and estimating the expansion in line \eqref{ineq:PDGab_in_CP_12_2} with the help of \eqref{ineq:PDGab_pos}
leads to the 
inequality:
\begin{align*} 
\triangle_k \left(\xbar, \pbar \right) 
\geq &\, \left(\gamma + \frac{1}{\tau_k} \right)  \BregD_{\X}\left(\xbar, x_{k+1} \right)  +
\frac{\BregD_{\Y^*}\left(\pbar, p_{k+1} \right)}{ \sigma_k} 
+ \frac{\BregD_{\Y^*} \left( p_{k+1}, p_{k}  \right)}{ \sigma_k} 
 + \frac{\BregD_{\X} \left(x_{k+1}, x_{k} \right)}{ \tau_k} \\
\ & + \left< \operator (x_{k+1} - \xhat_k ), p_{k+1} - \pbar \right>_{\Y} .
\end{align*}
Now we use Lemma \ref{cor:Young_BS} with \( \alpha = \frac{C_{\X}}{\|\operator\|\, \sigma_k\, \tau_k} \)  
\[
- \theta_{k-1}\, \left< \operator \left(x_{k} - x_{k-1} \right), p_{k+1} - p_k \right>_{\Y}  
\geq -\frac{\BregD_{\Y^*} \left( p_{k+1}, p_{k}  \right)}{ \sigma_k} 
-\frac{\theta_{k-1}^2 \, \| \operator\|^{2}\,  \tau_{k-1} \, \sigma_k}{C_{X}\, C_{Y^*}}
\frac{\BregD_{\X} \left(x_{k}, x_{k-1} \right)}{\tau_{k-1}}, 
\]
and insert \( \xhat_k = x_k + \theta_k (x_k - x_{k-1}) \) from \eqref{alg:CP_BS_3} such that we end up with
\begin{align*} 
\triangle_k \left(\xbar, \pbar \right) 
 \geq & \left( 1 + \gamma\, \tau_k  \right)\frac{\tau_{k+1}}{\tau_k} \frac{\BregD_{\X} \left( \xbar, x_{k+1} \right)}{ \tau_{k+1}}
 + \frac{\sigma_{k+1}}{\sigma_k} \frac{\BregD_{\Y^*}\left(\pbar, p_{k+1} \right)}{ \sigma_{k+1}} \\
 &+ \frac{\BregD_{\X} \left(x_{k+1}, x_{k} \right)}{ \tau_k}
- \frac{\theta_{k-1}^2 \, \| \operator\|^{2}\,  \tau_{k-1} \, \sigma_k}{C_{X}\, C_{Y^*}} 
\frac{\BregD_{\X} \left(x_{k}, x_{k-1} \right)}{\tau_{k-1}} \\
\ & + \left< \operator (x_{k+1} - x_k ), p_{k+1} - \pbar \right>_Y 
 - \theta_{k-1} \left< \operator(x_{k} - x_{k-1} ), p_{k} - \pbar \right>_Y . 
\end{align*}
The choice of the parameters ensures that
\[ 
\left( 1 + \gamma\, \tau_k \right)\frac{\tau_{k+1}}{\tau_k} 
= \theta_k^{-1}\geq 1 , \quad 
\frac{\sigma_{k+1}}{\sigma_k} 
= \theta_k^{-1} \geq 1, 
\quad \text{and} \quad 
\frac{\tau_{k}}{\tau_{k+1}} 
=  \theta_k^{-1} \geq 1 
\quad \text{for all } k \in \Nset.
 \]
Moreover, because of 
\( \min \left( C_X, C_{Y^*} \right) \geq \|\operator \|\ \tau_0^{\frac{1}{2\q}} \sigma_0^{\frac{1}{2\qstar}} 
= \|\operator\|\ \tau_k^{\frac{1}{2\q}} \sigma_k^{\frac{1}{2\qstar}}\) 
for \( k \in \Nset \) we have
\begin{align*} 
\frac{1}{\tau_{k}} \frac{\theta_{k-1}^2 \, \| \operator\|^{2}\,  \tau_{k-1} \, \sigma_k}{C_{X}\, C_{Y^*} } 
= \frac{1}{\tau_{k-1}} \frac{ \| \operator\|^{2}\,  \tau_{k} \, \sigma_{k}}{C_{X}\, C_{Y^*}} \leq \frac{1}{\tau_{k-1}}.
\end{align*}
Therefore 
\begin{align*} 
\frac{\triangle_k \left(\xbar, \pbar \right) }{\tau_k}
\geq &  \frac{\triangle_{k+1}\left(\xbar, \pbar \right) }{\tau_{k+1}}
 + \frac{\BregD_{\X} \left(x_{k+1}, x_{k} \right)}{ \tau_k^2} 
- \frac{\BregD_{\X} \left(x_{k}, x_{k-1} \right)}{\tau_{k-1}^2} \\
&+\frac{1}{\tau_k} \left< \operator(x_{k+1} - x_k ), p_{k+1} - \pbar \right>_Y 
 - \frac{1}{\tau_{k-1}}  \left< \operator (x_{k} - x_{k-1} ), p_{k} - \pbar \right>_Y 
\end{align*}
holds. Now, summing these inequalities from \( k = 0 \) to \(N-1\) for some \( N > 0 \) with \( x_{-1} \coloneqq x_0 \) 
and applying \eqref{ineq:apply_Young_ineq_sum} with \( \tau = \tau_{N-1} \)
leads to
\begin{align*}
\frac{\triangle_0 \left(\xbar, \pbar \right)}{\tau_0}
&\geq  \frac{\triangle_{N}\left(\xbar, \pbar \right) }{\tau_{N}}
 + \frac{\BregD_{\X} \left(x_{N}, x_{N-1} \right)}{ \tau_{N-1}^{2 \qstar}}  
+\frac{1}{\tau_{N-1}} \left< \operator(x_{N} - x_{N-1} ), p_{N} - \pbar \right>_Y \\
&\geq \frac{\triangle_{N}\left(\xbar, \pbar \right) }{\tau_{N}} 
- \frac{\|\operator\|^{2 \qstar}}{C_{\X} \, C_{\Y^*}} \BregD_{\Y^*} \left( \pbar,  p_{N} \right).
\end{align*}
By multiplying by \( \tau_N^{2} \)  and using the identity 
\( \tau_N\, \sigma_{N} = \tau_0\, \sigma_{0} \)
we obtain the following error bound:  
\[ \frac{\tau_{N}^{2}}{\tau_0\, \sigma_{0}} \left( 1- \frac{||\operator ||^{2 \qstar}}{ C_{\X} \, C_{\Y^*}} \tau_0\, \sigma_{0} \right)
\BregD_{\Y^*}\left(\pbar, p_{N} \right) + \BregD_{\X} \left( \xbar, x_{N} \right)
\leq \tau_N^{2\qstar} \left(\frac{\BregD_{\Y^*}\left(\pbar, p_{0} \right)}{ \sigma_0\,\tau_0 } 
+ \frac{\BregD_{\X} \left( \xbar, x_{0} \right)}{ \tau_0^{2}}\right).\]
Substituting \( \gamma \) by \( \frac{\gamma}{2} \) in
Lemma 1-2 and Corollary 1 in \cite{ChambollePock} shows that for any \( \epsilon > 0 \) 
there exists a \( N_0 \in \Nset\) (depending on \(\epsilon\) and \( \gamma\, \tau_0 \)) 
with \( \tau_N^2 \leq 4(1+ \epsilon) (N\, \gamma)^{-2} \)
for all \( N \geq N_0 \). 
This completes the proof.
\end{myproof}

Note that compared to the error estimate in \cite[Theorem 2]{ChambollePock}
the error bound \eqref{eq:error_bound_CPBS2} is \( 4 \) times larger for the generalized version \( \CPBS \). 
That is due to the fact that in the Hilbert space case also the positive term \eqref{ineq:PDGab_pos}
is bounded from below by \(\frac{\gamma}{2} \| x_{k+1} -\xbar \|_{\X}^2 \). 
In the considered Banach space setting we obtain \( \gamma\, \BregD_{\X}( x_{k+1},\xbar )\) 
as a lower bound, while \( \gamma\, \BregD_{\X}(\xbar, x_{k+1})\) would be required in order to prove the same result.

Finally, we will show that under the additional assumption that \( g^* \) 
fulfills \eqref{ineq:Def_gereneralization_uniformlyCVX} for some \( \delta \) we will achieve linear convergence:
\begin{mythm}\label{thm:CP_3}
Assume both \(f \) and \( g^* \) to satisfy property \eqref{ineq:Def_gereneralization_uniformlyCVX} 
for some constants \( \gamma > 0 \)  and \( \delta > 0 \), respectively.
Then for a constant parameter choice 
\(  \sigma_k = \sigma, \tau_k  = \tau,\) \( \theta_k = \theta,  {k\in \Nset} \),
with
\begin{itemize}
  \item  \( \mu \leq   \frac{\sqrt{\gamma\, \delta}\min \left\{C_X, C_{Y^*} \right\}}{ \|\operator \|} ,\)
 \item \( \sigma = \frac{\mu}{ \delta},\quad \tau = \frac{ \mu}{ \gamma} ,  \)
 \item \(\theta \in \left[ \frac{1}{1+\mu},1 \right], \)
\end{itemize}
the sequence \( \left(x_k, p_k \right)_{k \in \Nset} \) 
we receive from algorithm \ref{alg:CP_BS1} has the error bound:
\begin{equation} 
 \left( 1- \omega \right) \delta\, \BregD_{\Y^*}\left(\pbar, p_{N} \right) 
 +\gamma\, \BregD_{\X} \left( \xbar, x_{N} \right)
\leq \omega^N \left( \delta \, \BregD_{\Y^*}\left(\pbar, p_{0} \right)
+  \gamma \, \BregD_{\X} \left( \xbar, x_{0} \right) \right),
\end{equation}
with \( \omega = \frac{1+\theta}{2 +\mu} < 1 \). 
\end{mythm}

\begin{myproof}
In analogy to the proof of Theorem \ref{thm:convergence_CP_BS_2}, 
 we obtain from property \eqref{ineq:Def_gereneralization_uniformlyCVX} of \( f \) and \( g^*\)  
 a sharper estimate for \eqref{ineq:CP_12_weighted_error}-\eqref{ineq:CP_12_r}, where we set 
 \( \left(x,y^* \right) = \left(\xbar, \ybar \right) \):
We replace \eqref{eq:def_subd_for_x_k+1} by \eqref{ineq:subd_fmidcvx} 
and \eqref{eq:def_subd_for_p_k+1} by 
\[ 
g^* \left( \pbar \right) - g^* \left( p_{k+1} \right)  \geq  
 \left< \frac{ J_{Y^*}\left(p_k \right)- J_{Y^*} \left(p_{k+1} \right)}{\sigma}, \pbar - p_{k+1} \right>_Y + 
\left< \operator \xhat_k, \pbar -p_{k+1} \right>_{Y}
+\delta\, \BregD_{\Y^*}\left(\pbar, p_{k+1} \right).
\]
This together with \eqref{ineq:PDGab_pos}
leads to
\begin{align} \label{ineq:CP_51}
\begin{split}
\triangle_k \left(\xbar, \pbar \right) 
\geq &  \left( \delta + \frac{1}{\sigma} \right) \BregD_{\Y^*} \left(\pbar, p_{k+1} \right)
+ \left(\gamma + \frac{1}{\tau} \right)  \BregD_{\X}\left(\xbar, x_{k+1} \right)
+ \frac{\BregD_{\Y^*} \left( p_{k+1}, p_{k}  \right)}{ \sigma} 
 + \frac{\BregD_{\X} \left(x_{k+1}, x_{k} \right)}{ \tau} \\
\ & + \left< \operator (x_{k+1} - \xhat_{k} ), p_{k+1} - \pbar \right>_{\Y} .
\end{split}
\end{align}
Using \eqref{alg:CP_BS_3} and \eqref{ineq:Young_BS} with some \( \alpha > 0 \) 
we can estimate the last term in the following way:
\begin{align*}
 \left< \operator (x_{k+1} - \xhat_{k} ), p_{k+1} - \pbar \right>_{\Y}
 =& \left< \operator \left(x_{k+1} - x_k \right), p_{k+1} - \pbar \right>_{\Y}
 - \omega \left< \operator \left(x_{k} - x_{k-1} \right), p_{k} - \pbar \right>_{\Y} \\
  &- \omega \left< \operator \left(x_{k} - x_{k-1} \right), p_{k+1} - p_k \right>_{\Y}
 - \left(\theta - \omega \right) \left< \operator \left(x_{k} - x_{k-1} \right), p_{k+1} - \pbar \right>_{\Y} \\
 \geq & \left< \operator \left(x_{k+1} - x_k \right), p_{k+1} - \pbar \right>_{\Y}
 - \omega \left< \operator \left(x_{k} - x_{k-1} \right), p_{k} - \pbar \right>_{\Y} \\
  &- \omega\, \|\operator \| \frac{\BregD_{\Y^*}\left(p_{k+1}, p_{k} \right)}{C_{\Y^*}\, \alpha}
 - \theta\, \|T\|\, \alpha \frac{\BregD_{\X}\left(x_k, x_{k-1} \right)}{C_{\X}} \\ 
  &- \left(\theta - \omega \right)\|T\| \frac{\BregD_{\Y^*}\left( \pbar, p_{k+1} \right)}{C_{\Y^*}\, \alpha},
\end{align*}
for any \( \omega \in [(1+\mu)^{-1}, \theta] \).
Now, we set \( \alpha = \omega \left(\frac{\gamma}{\delta} \right)^{\frac{1}{2}}\)  such that 
\( \frac{\|T\|\, \mu\, \omega }{C_{\Y^*} \, \alpha} \leq \delta = \frac{\mu}{\sigma}\) 
and \(\frac{\mu\, \|T\|\, \alpha}{C_{\X} } \leq \omega\, \gamma \)
and multiply inequality \eqref{ineq:CP_51} by \( \mu\):
\begin{align} 
\label{ineq:CP_57}
\mu \, \triangle_k \left(\xbar, \pbar \right) 
\geq &  \left( 1 + \mu - \frac{1}{\omega} \right)\, \mu \, \triangle_{k+1} \left(\xbar, \pbar \right) 
+  \frac{\mu}{\omega} \triangle_{k+1} \left(\xbar, \pbar \right) 
 + \gamma\, \BregD_{\X} \left(x_{k+1}, x_{k} \right)
 -  \theta\, \omega\, \gamma\, \BregD_{\X} \left(x_{k}, x_{k-1} \right) \\
\ & + \mu\, \left< \operator \left(x_{k+1} - x_k \right), p_{k+1} - \pbar \right>_{\Y}
 - \mu\, \omega \left< \operator \left(x_{k} - x_{k-1} \right), p_{k} - \pbar \right>_{\Y}
  - \frac{\left(\theta - \omega\right) \, \delta}{ \omega} \BregD_{\Y^*}\left( \pbar, p_{k+1} \right).
\nonumber
 \end{align}
As in \cite{ChambollePock} we choose 
\begin{equation}\label{def:omega}
\omega = \frac{1+\theta}{2+\mu}  
\geq \frac{1+\theta}{2+\frac{\sqrt{\gamma\, \delta}\min \left\{C_X, C_{Y^*} \right\}}{ \|\operator \|}},   
\end{equation}
in order to ensure that
\[ 
\left( 1 + \mu - \frac{1}{\omega} \right)\, \mu \, \triangle_{k+1} \left(\xbar, \pbar \right)
 - \frac{\left(\theta - \omega\right) \, \delta}{ \omega} \BregD_{\Y^*}\left( \pbar, p_{k+1} \right)
 \geq 0 .
 \]
Thus, multiplying \eqref{ineq:CP_57}
with \( \omega^{-k} \) and summing from \(k= 0 \) to \( N-1\) for some \( N > 0 \) 
 where we set \( x_{-1} = x_0 \)) leads to 
\begin{align*} 
\mu \, \triangle_0 \left(\xbar, \pbar \right) 
\geq  \omega^{-N}  \mu \, \triangle_{N} \left(\xbar, \pbar \right) 
 + \omega^{-N+1}\, \gamma \BregD_{\X} \left(x_{N}, x_{N-1} \right) 
+\omega^{-N+1}\, \mu\, \left< \operator \left(x_{N} - x_{N-1} \right), p_{N} - \pbar \right>_{\Y}.
\end{align*}
Finally, by using Lemma \ref{cor:Young_BS} with 
\( \alpha = \left(\nicefrac{\gamma}{\delta}\right)^{\frac{1}{2}} \), 
we obtain from 
\( \| T \|\, \mu\, \alpha \leq \gamma\, \min \left\{C_{\X}, C_{\Y^*} \right\} \)
as well as 
\( \nicefrac{\| T \|\, \mu}{ \alpha} \leq \delta\, \min \left\{C_{\X}, C_{\Y^*} \right\}\):
\begin{align*} 
\mu \, \triangle_0 \left(\xbar, \pbar \right) 
\geq  \omega^{-N}  \mu \, \triangle_{N} \left(\xbar, \pbar \right) 
 - \omega^{-N+1}\, \delta\, \BregD_{\Y^*} \left(\pbar,  p_{N} \right),
\end{align*}
which completes the proof.
\end{myproof}

\begin{myremark}
 Because of \eqref{ineq:Bregman_pNorm} we proved in Theorem \ref{thm:convergence_CP_BS_2}
 a convergence rate of \(\mathcal{O} \left( \frac{1}{N}\right) \), 
 while Theorem \ref{thm:CP_3} even gives a convergence rate of 
 \(\mathcal{O} \left( \omega^{\frac{N}{2}} \right) \), 
 i.e. linear convergence.
\end{myremark}

\begin{myremark} \label{remark_relaxing_C_XY}
The parameter choice rules provided by Theorems \ref{thm:CP_BS1}, \ref{thm:convergence_CP_BS_2} 
and \ref{thm:CP_3} depend on the constants \( C_{\X} \) and \( C_{\Y^*} \) 
given by \eqref{ineq:Bregman_pNorm}. 
That is due to the application of Lemma \ref{cor:Young_BS} in the 
corresponding proofs. Now, if we assume that 
for a specific application this estimate 
is only required on bounded domains, the constants  
\( C_{\X} \) and \( C_{\Y^*} \) might be not optimal.
In fact, the numerical experiments indicate that we
obtain faster convergence if we relax the parameter choice of 
\( \sigma\) and \( \tau \) 
by replacing the product  \( C_{\X}\, C_{\Y^*} \) by a value
\(C \in \left[ C_{\X}\, C_{\Y^*} , 1  \right]\) close to 1.
\end{myremark}

\section{Duality mappings and generalized resolvents} \label{sec:general_resolvents}

In this section we give some examples of duality mappings and discuss the special generalization of the resolvent in our algorithm. 
\paragraph*{Duality mappings}
As shown in \cite{Hanner}, the reflexive Banach space \( l^r \) with \( r \in \left(1, \infty \right) \) 
is \(\max \left\{r,2\right\}\)-convex and \( \min \left\{r, 2 \right\}\)-smooth. 
One easily checks that the same holds true for the weighted sequence space \( l^r_W \) 
with positive weight \( W \) and norm
\( \left\| x \right\|_{l^r_{W}} \coloneqq \left( \sum \limits_{j} w_{j} |x_{j}|^r \right)^{\frac{1}{r}} 
=  \left\|  \left(w_{j}^{\frac{1}{r}} x_{j} \right)_{j}\right\|_{l^r}\).
With respect to the  \(l^2\)-inner product the dual space is given by \( Z^* = \mathnormal{l}^{r^*}_{W^{-1}} \)
where \( W^{-1} = (w_j^{-1}) \)
and we have
 \[ J_{q,\mathnormal{l}^{r}_W}: \mathnormal{l}_W^r \to  \mathnormal{l}_{W^{-1}}^{r^*}, \quad J_{q,l^r_W} (x)  = \frac{W}{|| x||_{l^r_W}^{r-q}} |x|^{r-1}\ \text{sign}(x), \]
 which has to be understood componentwise.

In order to model for example ``blocky'' structured solutions \(\xbar \in  \X \)  
let us consider (a discretization of) a Sobolev space
$H^{1,r}\left( \torus^d \right)$ where \( \torus = (-\pi, \pi)\) 
with periodic boundary conditions or equivalently the unit circle $S^1$. 
To define these spaces we introduce 
the Bessel potential operators $\Lambda_s:=(I-\Delta)^{-s}$ by 
\[
\Lambda_s\phi := \sum_{n\in\mathbb{Z}^d}(1+|n|^2)^{-s/2}
\widehat{\phi}(n) \exp(in\cdot), \qquad s\in\mathbb{R},
\] 
a-priori for $\phi\in C^{\infty}(\torus^d)$ where 
$\widehat{\phi}(n):=(2\pi)^{-d}\int_{\torus^d}\exp(-inx)\phi(x)\,dx$ denote  
the Fourier coefficients. Note that $\Lambda_0=I$ and 
$\Lambda_s\Lambda_t =\Lambda_{s+t}$ for all $s,t\in\mathbb{R}$. 
For $s\geq 0$ and $r\in (1,\infty)$ the operators $\Lambda_{-s}$ 
have continuous 
extensions to $L^{r}(\torus^d)$, and so the Sobolev spaces 
\[
H^{s,r}(\torus^d):= \Lambda_{-s}L^r(\torus^d)
\qquad \mbox{with norms}\qquad \|\phi\|_{H^{s,r}(\torus^d)}:= \|\Lambda_s\phi\|_{L^r(\torus^d)}
\]
are well defined.  Actually, this definition also makes sense for $s<0$, 
and we have the duality relation 
\[ \left(H^{r,s}(\torus^d)\right)^* =  H^{-s, r^*}(\torus^d) 
\]
for $1/r+1/r^*=1$ (see e.g.\ \cite[\S 13.6]{taylor3:96}).
The normalized duality mapping 
\( J_{ H^{s,r}(\torus^d)}: H^{s,r}(\torus^d)  \to H^{-s, r^*}(\torus^d)\) 
is given by 
\[ J_{ H^{s,r}(\torus^d)} = \Lambda_{-s}\ J_{L^r}\ \Lambda_{s}.
\] 
Recall that for $s\in\mathbb{N}$ the space $H^{s,r}(\torus^d)$ coincide 
with the more commonly used Sobolev spaces 
\[ W^{s,r} \left( \torus^d \right) \coloneqq
 \left\{ \phi\in L^r(\overline{\torus^d}) \, | \   D^{\alpha} \phi \in L^r(\torus^d) \mbox{ for all }|\alpha|\leq s\right\}, 
\]
with equivalent norms (\cite[\S 13.6]{taylor3:96}). 
\( H^{s,r}\left( \torus^d \right)  \) is a separable, reflexive,  
\( \max \left\{r,2 \right\}\)-convex and \( \min \left\{r,2 \right\}\)-smooth Banach space (see e.g. \cite{adams2003sobolev}, \cite{XuRoach} ).

In the discrete setting we approximate $\torus$ by the grid
\(
\torus_N:=\frac{\pi}{N}\left\{-\frac{N}{2},-\frac{N}{2}+1,\dots,\frac{N}{2}-1\right\}
\)
for some $N\in 2\mathbb{N}$. The dual grid in Fourier space is 
$\widehat{\torus_N}:=\frac{N}{\pi}\torus_N$, and the discrete Fourier transform 
$\FT_N:=(2\pi N)^{-d}(\exp(ix\xi))_{x\in\torus_N^d,\xi\in\widehat{\torus_N}^d}$ and 
its inverse $\FT_N^{-1}$ can be implemented by FFT. Hence, the 
Bessel potentials are approximated by the matrices 
\[
\Lambda_{s,N}:= \FT_N^{-1}\operatorname{diag}
\left[(1+|\xi|^2)^{-\frac{s}{2}}:\xi\in\widehat{\torus_N}^d\right]\FT_N.
\]
On the finite dimensional space of grid functions 
$\phi_N:\torus_N^d\to \mathbb{C}$ we introduce the norms  
\[
\|\phi_N\|_{H^{s,r}(\torus_N^d)}:= \|\Lambda_{s,N}\phi_N\|_{l^r(\torus_N^d)}.
\] 
If $\torus_N^d$ is replaced by $a\torus_N^d$ for some $a>0$, 
then $\widehat{\torus_N^d}$ has to be replaced by 
$\frac{1}{a}\widehat{\torus_N}^d$.

\paragraph*{Generalized resolvents}
Setting \( F = J_{X}\) the resolvent \( \left( \partial f + J_{X}\right)^{-1} \) is obviously  
closely related to the \( F\)-resolvents \( (A+F)^{-1} F \) of maximal monotone operators \( A \)
as used in \cite{KamimuraTakahashi_PPA_BS} and studied in 
\cite{Bauschke10generalresolvents}. 
Our focus lies on the evaluation of these operators.
The following generalization of Moreau's decomposition 
(see e.g. \cite[Theorem 31.5]{Rockafellar_cvxAnalysis}) 
allows us to calculate the generalized resolvent 
\( \left( \sigma_k\, g^* + J_{Y^*} \right)^{-1} \) in line \eqref{alg:CP_BS_1}
without knowledge of \( g^* \). 
This identity can be also derived from \cite[Theorem 7.1]{Bauschke10generalresolvents},
but for the convenience of the reader we present a proof for this special case:  

\begin{lemma}\label{thm:general_Moraus_id}
For any \( g \in \Gamma(\Y) \), \( \sigma > 0 ,\) and \( y \in \Y \) the minimization problem
\(\text{argmin}_{z \in Y} \left(\frac{\sigma}{2} \| z - \frac{y}{\sigma}  \|_{\Y}^2 + g(z) \right)\)
has a unique solution \( \ybar\) which is equivalently characterized by 
\( J_{\Y}(y -\sigma\, \ybar) \in \partial g(\ybar) \).
Therefore, the operator \( \left( J_{\Y^*} \circ \partial g + \sigma  I \right)^{-1}: \Y \to \Y \) is well defined
and single valued.
Moreover, the following identity holds:
\begin{equation} \label{eq:general_Moraus_id} 
\left(  \sigma \partial  g^* +  J_{\Y^*} \right)^{-1} \left( y \right)
=  J_{\Y} \left(y - \sigma\, \left( J_{\Y^*} \circ \partial g + \sigma  I \right)^{-1} 
\left( y \right) \right), \quad y \in \Y 
\end{equation}
\end{lemma}
\begin{myproof}
The first assertion follows from \cite[Theorem 2.5.1]{Zalinescu_cvxAna}, \cite[Theorem 3.4]{Cioranescu}
and the optimality condition 
\mbox{\( 0 \in \partial \left( \frac{\sigma}{2} \| \ybar - \frac{y}{\sigma}  \|_{\Y} + g(\ybar) \right)\)}.
In order to prove the second assertion, 
we set 
\[ \ybar = \left( J_{\Y^*} \circ \partial g + \sigma  I \right)^{-1}(y)\] 
for some \( y \in \Y \).
Moreover, let  \( \pbar \in \Y^* \) be a solution to the minimization problem
\begin{align*} 
\pbar 
= \argmin \limits_{p \in \Y^*} \left(  \sigma g^*(p) - \left< y, p \right>_{\Y} + \Phi_{2, \Y^*}(p) \right) 
= \argmax \limits_{p \in \Y^*} \left(   -g^*(p) + \left< \frac{y}{\sigma}, p \right>_{\Y} -  \frac{1}{\sigma} \Phi_{2,\Y^*}(p) \right),
\end{align*}
with \( \Phi_{2,Y^*}(p) \coloneqq \frac{1}{2} \|p\|_{Y^*}^{2} \).
Then \( \pbar \) can be rewritten as \( \pbar =  \left( \partial  \sigma g^* +  J_{\Y^*} \right)^{-1} \left( y \right) \),
cf. \eqref{eq:Gresolvent_argmin}.
Because of 
\[ 
\left(\frac{1}{\sigma} \Phi_{2,\Y^*}\right)^*(y) 
= \left<y ,\sigma\, J_{\Y}(y) \right>_{\Y} - \frac{\sigma^{2}}{2\, \sigma} \| J_{\Y} (y) \|_Y^{2} 
= \sigma\, \Phi_{2,\Y} (y)
\]
\( \ybar = \text{argmin}_{z \in Y} \left(\sigma \Phi_{2, Y } \left( z - \frac{y}{\sigma}\right) + g( z) \right) \)
is the solution \( \ybar \in Y \) to the corresponding Fenchel dual problem (cf. \eqref{prbl:primal}, \eqref{prbl:dual}). 
Now, \eqref{eq:extremal_relation} implies 
\( -\pbar \in \sigma\, \partial  \Phi_{2, \Y}\left( \ybar - \frac{y}{\sigma}\right) = \sigma J_{\Y} \left(\ybar - \frac{y}{\sigma}\right) \).
Thus we end up with
\[ 
-\left( \partial  \sigma g^* +  J_{\Y^*} \right)^{-1} \left( y \right) = - \pbar
= \sigma J_{\Y} \left(\left( J_{\Y^*} \circ \partial g + \sigma  I \right)^{-1} \left( y \right) - \frac{y}{\sigma}\right) .
 \qedhere \]
\end{myproof}

Our algorithm appears to be predestined for the case that 
\( f\) and \( g \) are given by Banach space norms 
\[ f(x) \coloneqq \frac{1}{2} \|x \|_{Z}^2, \quad g(y) =  \frac{1}{2} \|y -y_0 \|_{W}^2\]
for reflexive, smooth and 2-convex Banach spaces \( Z \) and \( W^* \)
and some \( y_0 \in W \). The natural choice of the space $X$ and $Y$ 
in this case is \( X = Z \), \( Y = W \). Then, due to the theorem of Asplund 
and the generalization of Moreau's decomposition \eqref{eq:general_Moraus_id}
 the generalized resolvents of  \( f \) 
and \( g^*(y^*) = \frac{1}{2}\left( \| y^* \|_{\Y^*}^2 +\left<y_0, y^* \right>_Y \right) \)
reduce to the corresponding duality mappings:
\begin{align} 
\label{eq:Gresolvent_normXsqrt}
\left(\tau \partial f + J_{\X} \right)^{-1}(x^*)
&= \left(\tau J_{\X} + J_{\X} \right)^{-1}(x^*) = \frac{1}{\tau + 1} J_{\X^*}(x^*) \\
\label{eq:Gresolvent_normYstarsqrt}
\left(\sigma \partial g^* + J_{Y^*} \right)^{-1}(y) 
&=  J_{\Y} \left(y- \sigma \left( J_{\Y}\ J_{Y^*}(\cdot - y_0)+ \sigma I \right)^{-1}(y) \right) 
= J_{\Y} \left( \frac{y-\sigma\, y_0}{  \sigma+1} \right).
\end{align}
Moreover, as  we have
\begin{align*}
 f\left( u \right)- f(x) - \left<u-x, x^* \right>_X =  \BregD_{\X}(u, x), 
 \quad g^*\left( p \right)- g^*(y^*) - \left<y, p-y^*  \right>_{\Y} =  \BregD_{\Y}(y^*, p).
\end{align*}
for all \( u,x \in X ,\) \(x^* \in \partial f(x) = J_{\X}(x) \) 
and \(y^*,p \in \Y^*\), \(y = \partial J_{\Y^*}(p) \), the functions  
\( f \) and \( g^* \) satisfy
property \eqref{ineq:Def_gereneralization_uniformlyCVX} 
for all \( \gamma, \delta \in (0,1] \), respectively.

If \( \X \neq Z \) or \(\Y \neq W \), however, 
a system of nonlinear equations has to be solved 
in order to evaluate the resolvents in lines \eqref{alg:CP_BS_1} and \eqref{alg:CP_BS_2}. 
In general for all exponents \( r \in \left(1, \infty \right) \), these resolvents of 
\begin{equation}\label{eq:fg_conv_Ljpunov}
  f(x) \coloneqq \frac{1}{r}\| x \|_{\X}^r, \quad  g(y) \coloneqq \frac{1}{r} \| y -y_0 \|_{\Y}^r,\ y_0 \in \Y
\end{equation}
 become rather simple:
 \begin{mycor}
  For \( \sigma, \tau > 0 \) and  \(f, g\) given by \eqref{eq:fg_conv_Ljpunov}, we have
  \begin{align} 
  \left(\tau \partial f + J_{\X} \right)^{-1}(x^*) &= \frac{1}{\tau\, \alpha^{r-2}+1} J_{\X^*}(x^*), 
  \quad  x^* \in \X \\ 
  \left(\sigma \partial g^* + J_{Y^*} \right)^{-1}(y) &= \frac{1}{ \beta^{r-2}+ \sigma} J_{\Y} \left( y-\sigma\, y_0  \right),
  \quad y \in \Y, 
   \end{align}
   where \( \alpha \geq 0 \) is the maximal 
   solution of \( \tau\, \alpha^{r-1} + \alpha = \| x^* \|_{X^*} \) and 
   \( \beta \geq 0 \) the maximal solution  of \( \beta^{r-1} +\sigma \,\beta = \| y - \sigma y_0 \|_Y\).
\end{mycor}
\begin{myproof}
Setting \( x = \left(\tau \partial f + J_{\X} \right)^{-1}(x^*) \),
the identity  \( \partial f (x)= J_{r,\X}(x) = \| x \|_{X}^{r-2}\ J_{\X}(x) \) implies that
\( \left( \tau\ \|x\|_{\X}^{r-2} + 1 \right) J_{\X}(x) = x^*\), and thus
 \( \alpha = \|x\|_{\X} = \| J_{\X}(x) \|_{\X^*} \geq 0 \). Inserting \( \alpha \) proves the first assertion. 
 For the second one we set \( \tilde{y}  \coloneqq \left(J_{Y^*} \circ \partial g + \sigma I \right)^{-1} (y) \).
 Because of \(y - \sigma y_0 =  J_{\Y^*} \circ \partial g (\tilde{y}) + \sigma (\tilde{y} -y_0) 
 = \left(\| \tilde{y} - y_0 \|_{\Y}^{r-2} + \sigma \right) \left(\tilde{y} - y_0\right)\)
 we have \( \beta = \| \tilde{y} - y_0 \|_{\Y} \) and \( \tilde{y} = \frac{y + \beta^{r-2}\, y_0 }{\beta^{r-2}+ \sigma}\).
 Now,  by Lemma \ref{thm:general_Moraus_id} follows:
\[ \left(\sigma \partial g^* + J_{Y^*} \right)^{-1}(y) 
=  J_{\Y} \left(y- \sigma   \tilde{y} \right) 
= J_{\Y} \left( \frac{y-\sigma\, y_0}{  \beta^{r-2}+ \sigma} \right) \qedhere.\]
\end{myproof}

Also other standard choices of \( f \) for which the resolvent has a closed form, 
provide a rather simple form for \( \left(\tau \partial f + J_{\X} \right)^{-1} \) as well.

\begin{myexp} \label{expl:Gresolvent_indicator_l1}
Consider the indicator function \( \chi_C \in \Gamma(\X) \) of a closed 
convex set \( C \subset \X \):
 \[\chi_C(x) = \begin{cases} 
                0 \quad &x \in C \\
                +\infty \quad &\text{otherwise}.
               \end{cases} \]
Then we obtain for any \( x^* \in X^* \) and any positive \( \tau \)  
\begin{align*} 
\left( \tau\ \partial \chi_C +J_{\X} \right)^{-1} (x^*) 
 &= \argmin \limits_{z \in \X}  \left(\tau\, \chi_C(z) - \left<z, x^* \right>_{\X}  + \frac{1}{2} \|z \|_{\X}^2 \right) \\
 &= \argmin \limits_{z \in C}  \left(\frac{1}{2} \left\| J_{\X^*}(x^*)\right\|_{\X}^2 
 - \left<z, J_{\X}\left(J_{\X^*}(x^*) \right) \right>_{\X}  + \frac{1}{2} \|z \|_{\X}^2 \right)
 = \pi_{C} (x^*) 
 \end{align*}
 where \( \pi_C: X^* \to C\) with \( \pi_{C}(x^*) \coloneqq \argmin \limits_{z \in C} \BregD_{\X}(z,J_{\X^*}(x^*)) \) 
 denotes the generalized projection 
 introduced by Alber \cite{Alber_IsraelSeminar}. 
For \( f(x) = \| x \|_{l^1} \) the subdifferential is given by
\begin{align*}
 \partial f (x) = \begin{cases}
                   \text{sign}(x) \quad &x \neq 0 \\
                   [-1,1]    \quad &\text{otherwise}.
                  \end{cases}
\end{align*}
Therefore we have for any \( x^* \in \X^* \) and any \( \tau > 0 \) 
\[ \left( \tau \partial f + J_{\X} \right)^{-1} (x^*) 
= J_{\X} \left( \max \left\{|x^*| - \tau , 0 \right\}\, \text{sign}(x^*) \right) .\]
\end{myexp}


\section{Numerical examples}\label{sec:numerical_exam}
In this section, we will test the performance of the generalized Chambolle-Pock method for linear and 
nonlinear inverse problems \( \operator x = y \), i.e.\ solving  \eqref{eq:Tikhonov_functional} or \eqref{eq:IRNM}.
In most examples, \( \X \) and \( \Y \) are weighted sequence spaces 
\( l_W^r \left( I \right) \) with \( r \in (1, \infty) \), countable or finite index sets \( I \), 
and positive weight \( W \), for which the required operator norm \(\| \operator \| \) 
is calculated by the power method of Boyd \cite{Boyd_PowerM_BS}.  
Also when  \( \X \) is the discrete Sobolev space \( H^{s,r}(\mathbb{\torus}^d_N) \)  this method can be applied, 
since  the operators \( \operator: \X \to  \Y = l_W^r \left( I \right)\) and 
\( A \coloneqq T \Lambda_{-s}: l^r(\torus^d_N) \to \Y \) have the same norms.
For all versions of the algorithm, we relax the parameter choice of 
\( \sigma \) and \( \tau \) according to Remark \ref{remark_relaxing_C_XY}. \\

\begin{figure}[t]
\begin{center}
  \includegraphics[width=\textwidth]{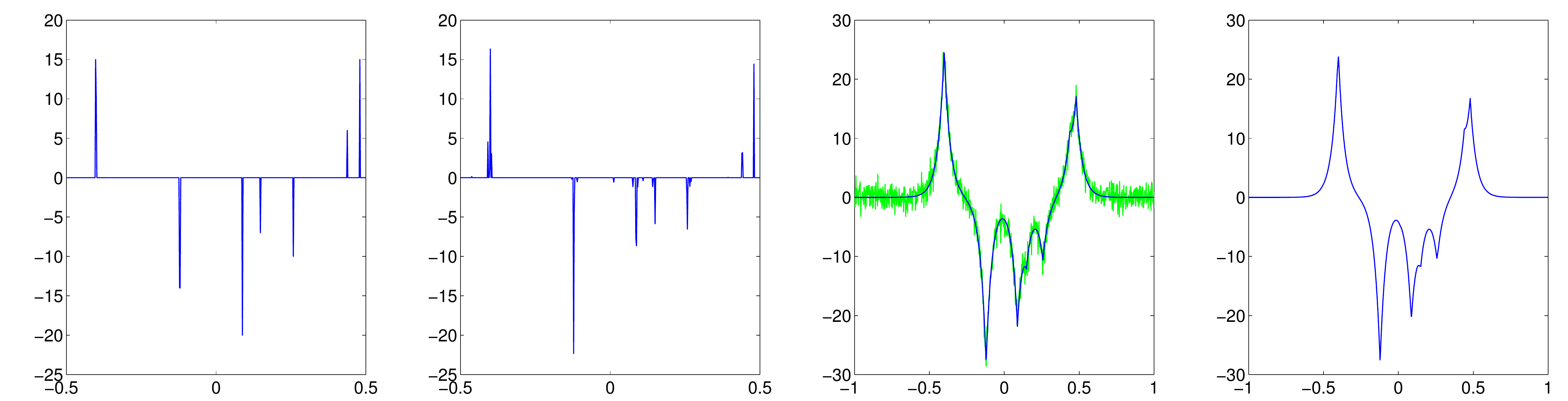}
 \caption{Deconvolution problem with penalty $R(x) = \|x\|_{l^1}$. From left to right: 
exact solution, reconstruction, exact (blue) and given (green) data, reconstructed data}
 \label{fig:conv_sol_data_recon}
 \end{center}
\end{figure}

\begin{figure}[t] 
 \includegraphics[width=1\textwidth]{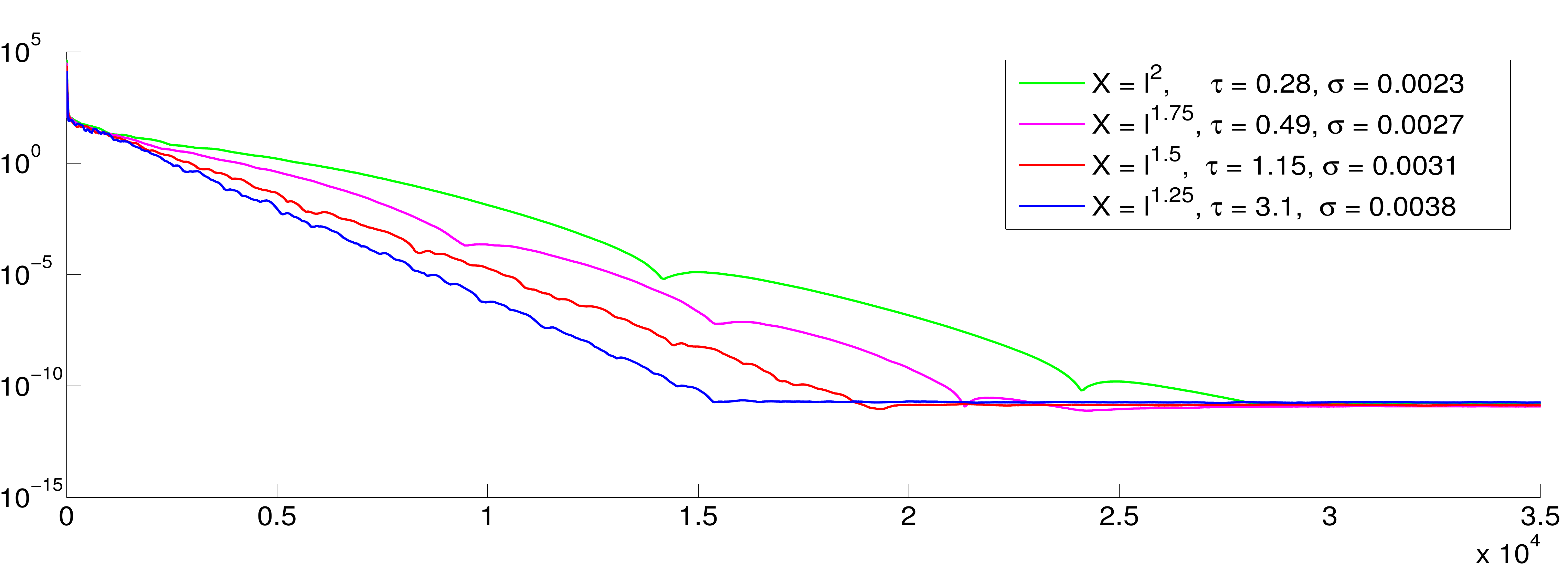}
 \caption{Convergence for the deconvolution problem with penalty $R(x) = \|x\|_{l^1}$. 
The error \( \|x_k - x_{\alpha} \|_{l^1}\) of the iterates \( (x_k)_{k \in \Nset} \) generated by the 
algorithm    \(\CPBS 1\) described in Theorem \ref{thm:CP_BS1} is plotted over the iteration step \( k \) 
for different choices of \(X \). The parameters \( \tau, \sigma \) are chosen optimally.}                                                                                                       
 \label{fig:convergence_conv}
\end{figure}

First, let us consider a linear ill-posed problem \( Tx = y,\) with convolution operator
\begin{equation} \label{eq:conv_operator}
 T(x):[-1,1] \to \Rset ,\quad T (x)(t) \coloneqq \int_{-\frac{1}{2}}^{\frac{1}{2}} x(s)\ k(t-s)\, ds, 
\quad k(t) \coloneqq \exp \left(- 5 |t| \right)
\end{equation}
and sparse solution \( x :\left[-\frac{1}{2},\frac{1}{2} \right] \to \Rset \) (see Figure \ref{fig:conv_sol_data_recon}).
This sparsity constraint is modeled by setting \( R(x) = \|x\|_{l^1} \) 
in \eqref{eq:Tikhonov_functional}. Moreover, as instead of the exact data \( y^{\dagger}\), 
only data \( y^{\delta} \) perturbed by 18 \% normal distributed noise is given,
we choose \( S(y^{\delta};y) =  \frac{1}{2}\| y^{\delta}-y \|_{l^2}^2 \) 
as data fidelity functional. According to the properties of the problem,
\( \X = l^r \left( I_{\X} \right), \) 
with \( r \in (1, 2] \) and \( \Y = l^2(I_{\Y}) \), seems to be a good choice. 
Here,
\(I_{\X} = \left\{-\frac{1}{2}, -\frac{1}{2}+\frac{1}{N-1},..., \frac{1}{2}  \right\} \)
is the discretization of  \( \left[\frac{1}{2}, \frac{1}{2} \right] \) 
and \(I_{\Y} = \left\{-1, -1+\frac{2}{N-1},..., 1 \right\} \) the of
\(\left[-1,1\right] \). The discretization of \(T\) is the discrete convolution. 
Now, for \(r = 2, 1.75, 1.5, 1.25 \) and \( \alpha = 5 \) 
we apply the version described in theorem \eqref{thm:CP_BS1} 
of our algorithm to \eqref{eq:Tikhonov_functional}. 
 Inspired by the optimality condition \( \operator \xbar \in \partial g^*(\pbar) \) 
with \( g^*(p) =\dataMisfit(y^{\delta}; p)^* =\frac{1}{2} \| p  \|_{l^2(I_{\Y})}^2 + \left<y^{\delta}, p \right>_{l^2}\),
 we pick
\begin{equation} \label{eq:inital_guess}
 \operator x_0 \in \partial g^*(p_0) = p_0 + y^{\delta} \Leftrightarrow p_0 =  \operator x_0 - y^{\delta}
\end{equation}
and \( x_0 = 0 \) as an initial guess.
The generalized resolvents are given by
 \eqref{eq:Gresolvent_normYstarsqrt} and example \ref{expl:Gresolvent_indicator_l1}.
Figure \ref{fig:convergence_conv} shows that for 
experimental optimal chosen parameters \( \sigma, \tau \) (according to Remark \ref{remark_relaxing_C_XY}),
we obtain faster convergence if  \( r\) turns 1. 
As Table \ref{tab:number_of_it} illustrates, this holds not only for the optimal parameter choice but also for any other choice of \( \sigma\). 
Here, we chose 
\(\tau \in (\sigma^{-1} \|T\|^{-2} -2^{-6}, \sigma^{-1} \|T\|^{-2} )\) for the Hilbert space case \( X = l^2(I_X), \)
and \( \tau \in \left[\sigma^{-1} \|T\|^{-2} C_1 ,\sigma^{-1} \|T\|^{-2} C_2 \right],\) with \(C_1 = 0.89,\ C_2 = 0.96 \in [0.25, 1], \)
for the Banach space case  \( X = l^{1.25} (I_X) \) (cf. Remark \ref{remark_relaxing_C_XY}).
Thus, we conclude that a choice of \( \X, \) which reflects the properties of the problem best, 
may provide the fastest convergence.

\begin{table}[tbt]
\centering
\begin{tabular}{|l|cl|cl|cl|}
\hline 
\multicolumn{1}{|c|}{ \(\mathbf{\sigma}\)} &  \( 0.007 \)& 	& \( 0.0023 \quad \)&	& \(  0.00075 \) & \\
\hline
\hline
\( \X = l^2(I_{\X})\) 		
& 76476&	(\(\tau = 0.0915\))	
& 22368& 	(\(\tau = 0.279\))	
& 39418&     	 (\(\tau = 0.854\)) \\
\(\X = l^{1.25}(I_{\X})\)	
& 26271&	(\(\tau = 1.7 \))	
& 16575&	(\(\tau = 4.8 \))	
& 38710&	(\(\tau = 14.66 \)) \\
\hline
\end{tabular}
\caption{Comparison of \( \CP \) with \( \X =l^2(I_{X} )\) and \(\CPBS 1 \) with \( X =l^{1.25}(I_X) \) 
for the deconvolution problem 
with different choices of \( \sigma \). The table shows the first iterations number $k$ 
for which \( \|x_{\alpha}-x_k \|_{l^1} \leq 10^{-5}\), averaged over 100 experiments.}
\label{tab:number_of_it}
\end{table}


\begin{figure}[th]
\begin{center}
\parbox{.55\textwidth}{
 \includegraphics[width=0.55\textwidth]{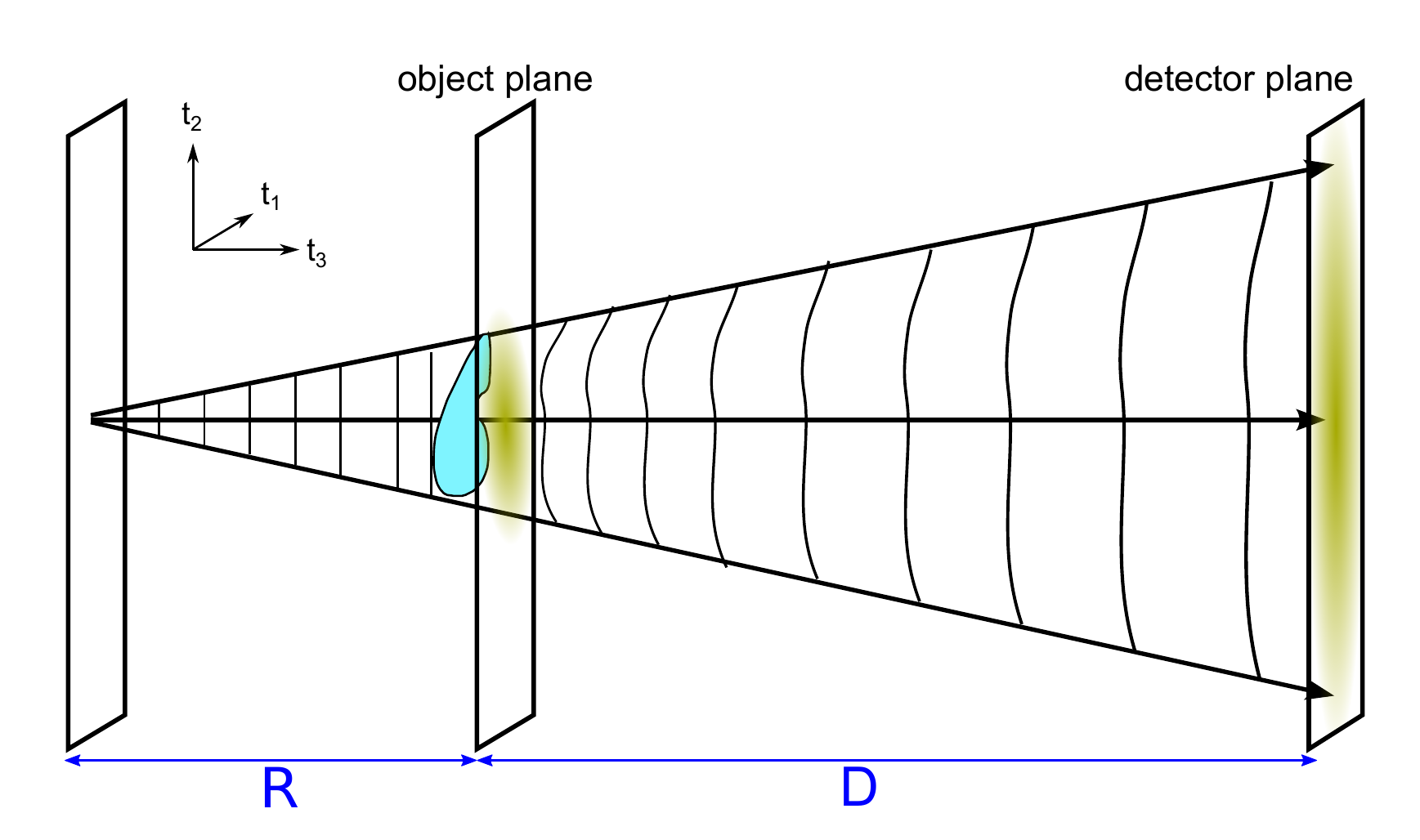}
 \caption{Experimental setup leading to the phase retrieval problem}
 \label{fig:phase_retrieval_setting}}
 \parbox{.4\textwidth}{\includegraphics[width=0.35\textwidth]{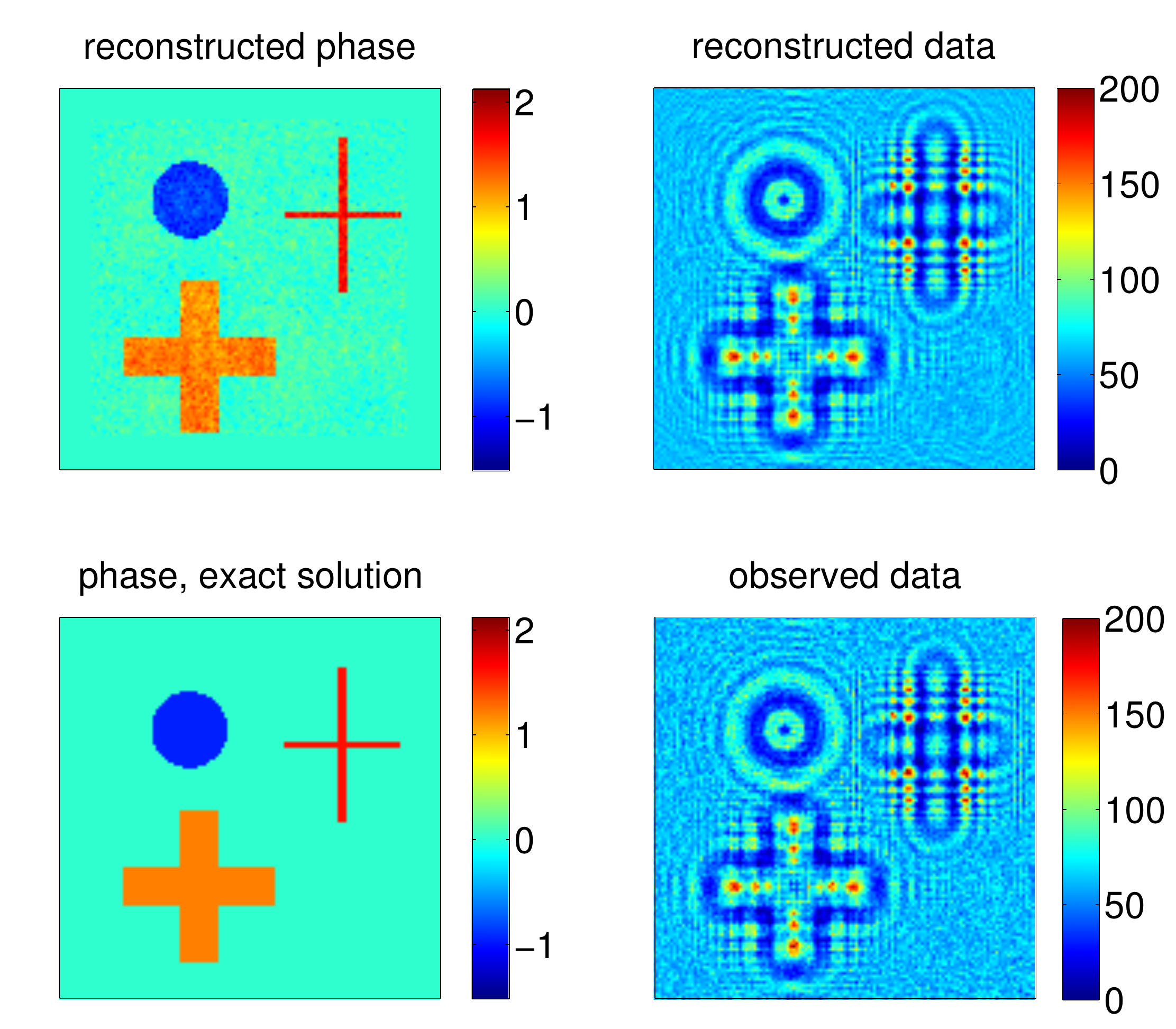}
 \caption{top: reconstructed phase $\phi$ and corresponding data $T(\phi)$ after 15 IRGN iterations;
bottom: exact phase $\phi^{\dagger}$ and simulated Poisson distributed diffracton pattern $y^{\delta}$ 
with $\mathbb{E}[y^{\delta}]= T(\phi^{\dagger})$}} 
 \end{center}
\end{figure}

\begin{figure}[th] 
 \includegraphics[width=1\textwidth]{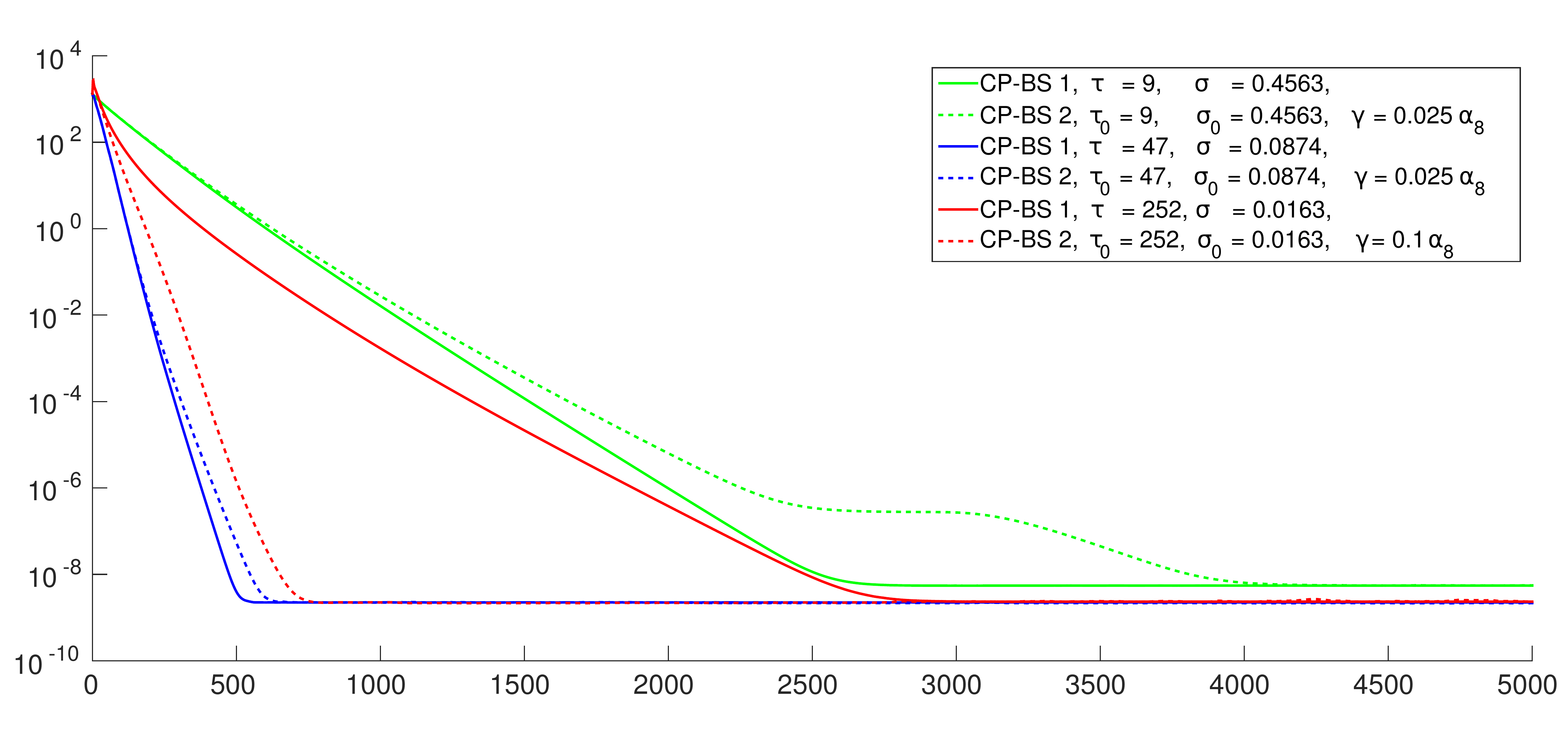}
 \caption{Convergence for the phase retrieval problem with penalty $R(x) = \frac{1}{2}\|x\|_{X}^2$ for $X=H^{1,1.1}(\torus^2)$ 
 at Newton step \( n=8\). The error \( \|x_n - x_k \|_{X}\) 
 of the iterates \( (x_k)_{k \in \Nset} \) of the algorithm \(\CPBS \) and the best approximation \( x_{n} \) to the true minimizer of 
 \eqref{eq:IRNM}  is plotted over the  iteration step \( k \).
 The parameter choice rules are defined by Theorem \ref{thm:CP_BS1} (\(\CPBS 1 \), solid) 
 and by Theorem \ref{thm:convergence_CP_BS_2} (\(\CPBS 2\), dotted), respectively.
 For given \( \tau \) (or \( \tau_0 \)), we set \( \sigma = 0.96 \| \operator^{\prime}[\phi_8]\|^{-2} \tau^{-1}\) (\( \sigma_0 \) analogously)} 
 \label{fig:convergence_phase_retrieval}
\end{figure}

As a second example, we consider a phase-retrieval problem (see figure \ref{fig:phase_retrieval_setting}): 
a sample of interest is illuminated by a coherent 
x-ray point source. From intensity measurements \( \left|u^{\delta}(\cdot, D) \right|^2 \)
of the electric field \( u: \Rset^3 \to \Cset \), which are taken in the detector plane,
orthogonal to the beam at a distance \( D > 0 \), 
we want to retrieve information on the refractive index of the sample.
More precisely, we are interested in the real phase 
\( \phi: \Rset^2 \times \left\{0 \right\} \to \Rset \) 
of the object function \( \obfct_{\phi}(x) = \exp \left(- \text{i} \kappa \phi(x) \right)\)
describing the sample, where \( \kappa \) denotes the wavenumber. 
We assume that $\e^{-\ima\, \kappa D}\, u(\cdot, D)$ 
can be approximated by the so called Fresnel propagator 
\( (P_{D}\obfct) \coloneqq  \FT^{-1}\left(\chi_{-\frac{D}{ \kappa}}
\cdot(\FT \obfct )\right)  \)
where \( \chi_{-c}(t_1,t_2) := \exp \left(- \text{i}\, c \left( t_1^2+t_2^2 \right) \right)\)  
is a chirp function with parameter \(c>0 \). 
Using the Fresnel scaling theorem we obtain the following forward operator $T$ mapping $\phi$ 
to  \( \left|u(\cdot, D) \right|^2 \): 
\[ (T(\phi))(M\, t_1,M\, t_2) 
\coloneqq \frac{1}{M^2} \left| \left(P_{\frac{D}{M}} \obfct_{\phi}\right) \left(t_1, t_2  \right) \right|^2 \,.
\label{eq:I2_eff}\]
Here \( M = \frac{R+D}{R} \), where \( R > 0 \) is the distance between the sample and the source,
denotes the geometrical magnification.
For a detailed introduction to this problem and phase retrieval problems in general we refer to \cite{Paganin}.
\( T \) is Fr\'{e}chet differentiable with
\[T^{\prime} [\phi](h)(M\, t_1,M\, t_2, D) 
= \frac{2}{M^2} \Re \left( \overline{P_{\frac{D}{M}}\left(\obfct(\phi)\right)(t_1,t_2)}\, 
\left(P_{\frac{D}{M}} \obfct^{\prime}_{\phi,h} \right) (t_1,t_2) \right) \]
so the IRNM \eqref{eq:IRNM} is applicable. This is an example with Poisson data,
thus following \cite{HohageWerner, WernerHohage}, after the discretization
\[ x = \left( x_j \right)_{j \in \torus^2_N} \coloneqq \phi(j_1, j_2)_{j \in \torus^2_N}, \quad
y = \left( y_j \right)_{j \in M\torus^2_N} \coloneqq |u(j_1, j_2, D)|^2_{j \in M\torus^2_N}
\]
we choose (using the convention \(0 \ln 0 \coloneqq 0 \)):
\begin{align*}  
S\left(y^{\delta}  ; y  \right)  =  
\begin{cases}
\sum \limits_{j \in M\torus^2_N} y_j - y_j^{\delta} \ln y_j  
\quad  
& \text{if } y\geq 0 \mbox{ and }y_j>0 \mbox{ for all } j\mbox{ with }y^{\delta}_j>0
\\
+\infty \quad &\text{otherwise}.
\end{cases} 
\end{align*}
Moreover, motivated by the weighted least square approximation (cf. \cite{Stueck12})
\[S\left( y^{\delta}  ; y  \right) 
\approx \frac{1}{2} \left\|\left( \frac{y_j^{\delta} - y_j}{\sqrt{y_j}} \right)_{j \in M\torus_N^2} \right\|_{l^2}^2 ,\]
in the \((n\!+\!1)\)-th iteration step of the IRNM, we consider the weighted space \( Y = l^2_{W} \) with weight
\( W \coloneqq \left(\operator(\phi_n) +\epsilon \right)^{-1}\)
and \( \epsilon =  0.1.\) Compared to setting \( Y =l^2 \), this leads to a faster convergence as numerical experiments show. 
\( J_{l^2_{W}} \left(y\right) = W\, y,  \) and \( \partial S\left(y^{\delta}  ; y \right)_j = 1 - \frac{y_j^{\delta}}{y_j} \) 
for \((y^{\delta} ,y) \in \text{dom } (S) = \left\{ (y^{\delta} ,y)\, |\ \dataMisfit(y^{\delta};y) < + \infty \right\}\)  imply
\[ \left(\left( J_{ l_W^{2}} \circ \partial S(y^{\delta}, \cdot) + \sigma I \right)^{-1} (y) \right)_j = \frac{y_{j}- W_{j}^{-1}}{2\, \sigma} 
+ \frac{\sqrt{\left(y_{j}-W_{j}^{-1} \right)^2 + 4\, \sigma\ W_{j}^{-1}\, y^{\delta}_j}}{ 2 \sigma}. \]
Hence, we obtain the generalized resolvent 
\( \left(\sigma\,\partial S(y^{\delta}, \cdot)^* + I \right)^{-1} \) 
by Lemma \ref{thm:general_Moraus_id}. 
The  "blocky" structured solution is taken into account by setting 
$X:= H^{1,r}(\torus^2)$ with $r=1.1$ and \( R(x) = \frac{1}{2} \| x\|_{X}^{2}.\) 
Note that although evaluating the generalized resolvent 
\( (\tau\, \alpha R + J_X)^{-1} = \frac{1}{\tau \, \alpha + 1 } J_{H^{-1, r^*}} 
= \frac{1}{\tau \, \alpha + 1 } \Lambda_{1} J_{l^{r^*}} \Lambda_{-1} \) 
is more expensive than in the case \( X = l^r \), it does not increase the complexity of the algorithm 
as the evaluations of \( T'[\phi] \) and \( T'[\phi]^* \) include Fourier transforms as well. 
Since \( R \) satisfies property \eqref{ineq:Def_gereneralization_uniformlyCVX},
we can apply the variant \(\CPBS 1 \) 
described in Theorem \ref{thm:CP_BS1} and also the variant \(\CPBS 2\) given by Theorem \ref{thm:convergence_CP_BS_2}.
Figure \ref{fig:convergence_phase_retrieval} compares both versions in the \( n = 8\)-th iteration step of the IRNM, where \( \alpha = 0.001 \).
The solid blue curve belongs to the version  \(\CPBS 1 \) for a optimal parameter choice of \( \tau \) and
\(\sigma \) we found experimentally. Note that for the limit \( \gamma \to 0 \) the parameter choice rule of \( \CPBS 2 \) 
coincide with the one of \( \CPBS 1 \). In fact, choosing \( \tau_0 \) and \( \sigma_0 \) in the same way as \( \tau \) and \( \sigma \)
that corresponds to this blue curve the version \( \CPBS 2 \) with \( \gamma = 0.0025 \alpha_8 \) gives the same curve. 
Tuning also the parameters \( \tau_0, \sigma_0  \) and \( \gamma \) (reasonable large)
in an optimal way, we did not obtain a better convergence result for \( \CPBS 2 \) than for \( \CPBS 1 \). 
However, \(\CPBS 2\) converges faster
for \( \tau = \tau_0 \) sufficiently large and adequately chosen \( \gamma \). 

\begin{figure}[t]
\begin{center}
 \includegraphics[width= 1\textwidth]{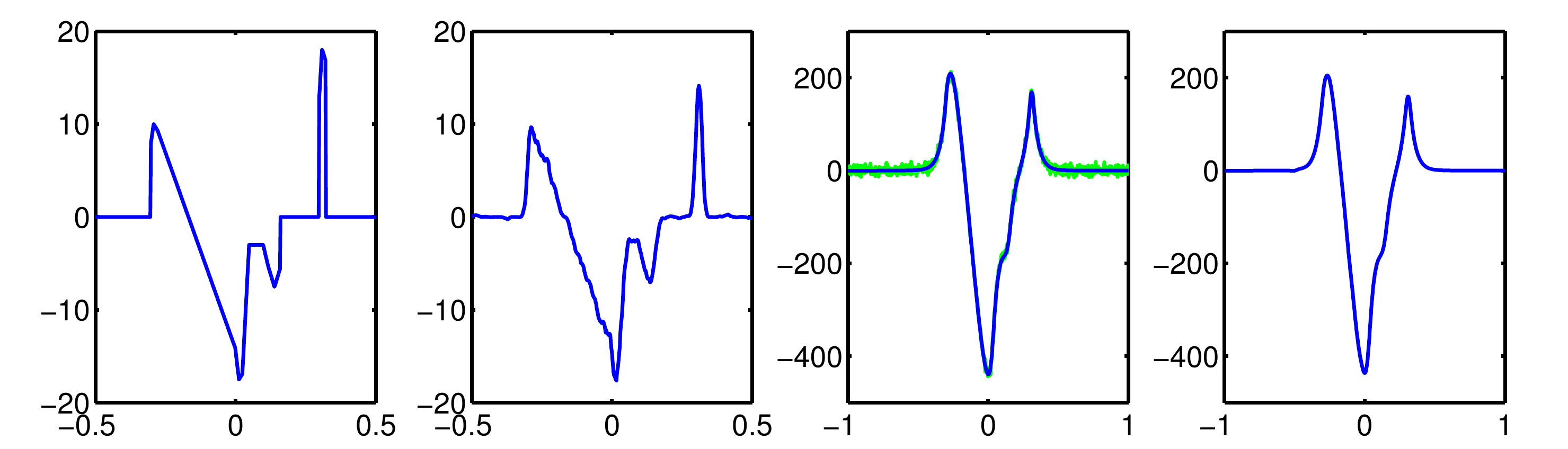}
 \caption{Deconvolution problem with penalty $R(x) =\frac{1}{2}\|x\|_{l^{1.5}}^2$. 
From left to right: 
exact solution, reconstruction, exact (blue) and given (green) data, reconstructed data}
 \label{fig:conv15_sol_data_recon}
 \end{center}
\end{figure}

In our last example, we apply the version described by Theorem \ref{thm:CP_3}, which we denote as \(\CPBS 3\), 
to the Tikhonov functional
\[ \frac{1}{2}\| \operator x - y^{\delta} \|_{\Y}^2 + \frac{\alpha}{2} \| x \|_{\X}^{2} \]
where \( T \) is again the convolution operator \eqref{eq:conv_operator}.  
We set \( \alpha = 1 \), \( X = l^{1.5}(I_X) \) and \( Y = l^2(I_X) \) (see figure \ref{fig:conv15_sol_data_recon}).
Setting \( \mu = \frac{C\, \sqrt{\gamma\, \delta}}{2\|T\|} \) with \( C = 0.98, \) we obtain for any choice
\( \gamma = \delta \in (0,2] \) the fastest convergence rate. The same rate is also provided by 
\( \CPBS 1 \) and \( \CPBS 2 \) for optimal chosen (initial) parameter. 
Compared to the first example where more than 15000 iterations were required to satisfy the stopping 
criterion, here we only need 558 iterations. 

\section*{Acknowledgement} We would like to thank Radu Bot and Russell Luke for helpful discussions. 
Financial support by DFG through CRC 755, project C2 is gratefully acknowledged. 

\bibliographystyle{abbrv}
\bibliography{./CP_BS}

\end{document}